\documentclass{amsart}

\usepackage{amsbsy,amssymb,amscd,amsfonts,latexsym,amstext,delarray,
amsmath,graphicx}
\input xypic
\usepackage{color}

\newtheorem{thm}{Theorem}[section]
\newtheorem{prop}[thm]{Proposition}

\newtheorem{lem}[thm]{Lemma}
\newtheorem{defn}[thm]{Definition}

\numberwithin{equation}{section}

\def\bR{{\mathbb R}}

\def\bT{{\mathbb T}}

\def\C{{\mathbb C}}
\def\F{{\mathbb F}}

\def\N{{\mathbb N}}

\def\Q{{\mathbb Q}}
\def\Z{{\mathbb Z}}
\def\R{{\mathbb R}}

\def\cA{{\mathcal A}}

\def\cQ{{\mathcal Q}}
\def\cR{{\mathcal R}}
\def\cS{{\mathcal S}}

\def\cW{{\mathcal W}}
\def\cX{{\mathcal X}}
\def\cY{{\mathcal Y}}
\def\cZ{{\mathcal Z}}

\title{Thermodynamic semirings}
\author{Matilde Marcolli and Ryan Thorngren}
\address{Mathematics Department, Mail Code 253-37, Caltech, 1200 E.~California Blvd. Pasadena, CA 91125, USA}
\email{matilde@caltech.edu} 
\email{rthorngr@caltech.edu}

\begin{document}
\maketitle

\begin{abstract}
The Witt construction describes a functor from the category of Rings to the category of characteristic 0 rings. It is uniquely determined by a few associativity constraints which do not depend on the types of the variables considered, in other words, by integer polynomials. This universality allowed Alain Connes and Caterina Consani to devise an analogue of the Witt ring for characteristic one, an attractive endeavour since we know very little about the arithmetic in this exotic characteristic and its corresponding field with one element. Interestingly, they found that in characteristic one, the Witt construction depends critically on the Shannon entropy. In the current work, we examine this surprising occurrence, defining a Witt operad for an arbitrary information measure and a corresponding algebra we call a thermodynamic semiring. This object exhibits algebraically many of the familiar properties of information measures, and we examine in particular the Tsallis and Renyi entropy functions and applications to nonextensive thermodynamics and multifractals. We find that the arithmetic of the thermodynamic semiring is exactly that of a certain guessing game played using the given information measure.
\end{abstract}

\tableofcontents

\section{Introduction}

The past few years have seen several interesting new results 
focusing on various aspects of the elusive ``geometry over the field with one element", 
see for instance \cite{Borg} \cite{CC} \cite{CCM} \cite{LP} \cite{Man}, \cite{Soule},
among many others. The idea of $\F_1$-geometry has its roots in an observation of
Tits \cite{Tits} that limits as $q\to 1$ of counting functions for certain varieties
defined over finite fields $\F_q$ exhibit an interesting combinatorial meaning,
suggesting that the resulting combinatorial geometry should be seen as an 
algebraic geometry over a non-existent ``field with one element" $\F_1$. 
Part of the motivation for developing a sufficiently refined theory of 
varieties and schemes over $\F_1$ lies in the idea that being able to 
cast ${\rm Spec} \Z$ in the role of a curve
over a suitably defined ${\rm Spec} \F_1$ may lead to finding an
analog for number fields of the Weil proof \cite{Weil} of the Riemann hypothesis 
for finite fields.

\smallskip

Among the existing approaches aimed at developing various aspects of
geometry over $\F_1$, the one that is of direct interest to us in the present paper
is a recent construction by Connes and Consani \cite{Co}, \cite{CC} of semirings
of characteristic one (a nilpotent hypothesis). These are endowed with an 
additive structure that provides an analog of the Witt formula for the addition of 
the multiplicative Teichm\"uller lifts in strict $p$-rings.  As observed in \cite{Co}
and \cite{CC}, the commutativity, identity, and associativity conditions for this addition 
force the function used in defining the Witt sums in characteristic one 
to be equal to the Shannon entropy.

\smallskip

The goal of this paper is to explore this occurrence of the Shannon entropy in 
the characteristic one Witt construction of \cite{Co} and \cite{CC}. 
In particular, we show here that the construction introduced in those papers
can be seen as part of a more general theory of ``thermodynamic semirings",
which encodes various properties of suitable ``entropy functions" in terms
of algebraic properties of the corresponding semirings. 

\smallskip

After reviewing the case of \cite{Co}, \cite{CC} in \S \ref{prelimSec}, 
we present a general definition and some basic properties of 
thermodynamic semirings in \S \ref{Shaxiomsec} and \S \ref{thermoSec},
based on the axiomatization of information-theoretic entropy through
the Khinchin axioms and other equivalent formulations. 
We then give in \S \ref{statmechSec} a physical interpretation of the
structure of thermodynamic semiring in terms of Statistical Mechanics,
distinguishing between the extensive and non-extensive cases and
the cases of ergodic and non-ergodic statistical systems. We see that
the lack of associativity of the thermodynamic semiring has
a natural physical interpretation in terms of mixing, chemical potentials,
and free energy. This generalizes the thermodynamic interpretation of 
certain formulas from tropical mathematics considered in \cite{Quad}.


\smallskip

We focus then on specific examples of other important information-theoretic entropy 
functions, such as the R\'enyi entropy, the Tsallis entropy, or the Kullback--Leibler
divergence, and we analyze in detail the properties of the corresponding
thermodynamic semirings.
In \S \ref{RenyiSec}, we consider the case of the R\'enyi entropy,
which is a one-parameter generalization of the Shannon entropy that still satisfies
the extensivity property. In \S \ref{tsallisSec}  we focus instead on the
Tsallis entropy, which is a non-extensive one-parameter generalization of
the Shannon entropy, and we show that a simple one-parameter deformation 
of the Witt construction of \cite{Co} and \cite{CC} identifies the Tsallis
entropy as the unique information measure that satisfies the associativity
constraint.

\smallskip

In \S \ref{KLsec} we consider the case of the Kullback--Leibler divergence or
relative entropy (information gain), and we show that thermodynamic
semirings based on this information measure can be associated to
univariate and multivariate binary statistical manifolds, in the sense of
information geometry, and to multifractal systems, in such a way that the
algebraic properties of the semirings detect the statistical
and multifractal properties of the underlying spaces. We also relate 
a hyperfield structure arising from the KL divergence to those considered in \cite{Vi}.
\smallskip

We also show in \S \ref{assocSec} that the algebraic structure of the thermodynamic semirings
can be encoded in a suitably defined successor function and that the properties
of this function and its iterates as a dynamical system capture both the
algebraic structure of the semiring and the thermodynamical properties of
the corresponding entropy measure. We give explicit examples of these
successor functions and their behavior for the Shannon, R\'enyi, and 
Tsallis entropies. In \S \ref{cumulsec} we show that this function has an 
interpretation as the cumulant generating function for the energy, which 
reveals some futher thermodynamic details of our construction.

\smallskip

Finally, in \S \ref{operadsec}, we phrase our construction 
using operads whose composition trees suggest an interpretation in terms of ``guessing games". 
Exploring this, we show that relations in a particular algebra--the thermodynamic semiring--for the guessing game operad 
correspond naturally to information theoretic properties of the entropy functions, cominiscent of 
an operadic characterization studied recently by Baez, Fritz and Leinster, which we review. 
This allows us to rephrase Connes and Consani's original construction in a way that makes clear 
why the Shannon entropy plays such a key role and provides a categorification of entropy functions.

\smallskip

In the last section we outline possible further directions, some of which
will eventually relate back the general theory of thermodynamic semirings
to the analogies between characteristic $p$ and characteristic one geometries.
Thus, this point of view based on thermodynamic semirings
may be regarded as yet another possible viewpoint on $\F_1$-geometry,
based on information theory and statistical geometry, a sort 
of ``cybernetic viewpoint".

\smallskip

\subsection{Witt vectors and their characteristic one analogs}

Witt vectors were first proposed by Ernst Witt in 1936 to describe unramified
extensions of the $p$-adic numbers. In particular, Witt developed integral
polynomial expressions for the arithmetic of strict $p$-rings in terms of
their residue rings. 

\smallskip

A ring $R$ is a strict $p$-ring when $R$ is complete and
Hausdorff under the $p$-adic metric, $p$ is not a zero-divisor in $R$, and the
residue ring $K = R / p R$ is perfect \cite{Lor}, \cite{Rab}, \cite{Serre}. 
The ring $R$ is determined by $K$
up to canonical isomorphism, and there is a unique multiplicative section $\tau : K
\rightarrow R$ of the residue morphism $\pi : R \rightarrow K$, ie.
\[ 
     \pi \circ \tau = {\rm id}_K, \ \ \ \  \tau (x y) = \tau (x) \tau (y) \ \ \ \   \forall x,
     y \in K.
   \]
Every element $x$ of $R$ can be written uniquely as
\[ 
x = \sum \tau (x_n) p^n, \ \ \ \  x_n \in K.
 \]
The $\tau (x)$ are called Teichm\"uller representatives. 

\smallskip

When $K
=\mathbb{F}_p$, $R =\mathbb{Z}_p$, but the Teichm\"uller representatives are
not $\{0, 1, \ldots, p - 1\}$ as they are in the common representation of
$\mathbb{Z}_p$. Instead they are the roots of $x^p - x$. We see from this
example that the arithmetic in terms of the Teichm\"uller representation above
is nontrivial. The Witt formula expresses the sum of these representatives as
\begin{equation}\label{TeichReps}
     \tau (x) + \tau (y) = \tilde{\tau} ( \sum_{\alpha \in I_p} w_p (\alpha,
    T) x^{\alpha} y^{1 - \alpha}),
\end{equation}
where $I_p =\{\alpha \in \mathbb{Q} \cap [0, 1]\, |\, p^n \alpha \in \mathbb{Z}$
for some $n\}$, $\tilde{\tau} : K [[T]] \rightarrow R$ is the unique map such that $\tilde{\tau} (x T^n) = \tau (x) p^n$, and $w_p
(\alpha, T) \in \mathbb{F}_p [[T]]$ is independent of $R$. Note that, since $K$
is perfect, the terms $x^{\alpha} y^{1 - \alpha}$ make sense.

\medskip

The idea of \cite{Co}, \cite{CC} is to generalize this to characteristic 
one by considering sums of the form
\begin{equation}\label{charonesums}
      x \oplus_w y := \sum_{\alpha \in I} w (\alpha) x^{\alpha} y^{1 - \alpha}
\end{equation}
where now $I =\mathbb{Q} \cap [0, 1]$ over sufficiently nice characteristic
one {\em semi}rings. 

\smallskip

According to Definition 2.7 of \cite{CC},
a semiring is characteristic one when $1 + 1 = 1$,
i.e. when it is idempotent. For example, the tropical semifield,
$\mathbb{T}=\mathbb{R} \cup \{ - \infty\}$, with addition given by the sup and
multiplication given by normal addition, forms a well studied characteristic
one semiring in the context of tropical geometry
\cite{ItMiShu}, \cite{MaStu}. 

\smallskip

Connes and Consani found in \cite{Co}, \cite{CC} that, over a suitably nice characteristic
one semiring, $\oplus_w$ is commutative, associative, shares an identity with
$+$, and is order preserving {\em if and only if} $w (\alpha)$ is of the form 
\begin{equation}\label{wSh}
w(\alpha) = \rho^{{\rm Sh} (\alpha)}, 
\end{equation}
where $\rho \geqslant 1$ and ${\rm Sh} (p)$ is the well known Shannon entropy
\begin{equation}\label{Shentropy}
  {\rm Sh} (p) = -C ( p \log p + (1 - p) 
\log (1 - p) ),
\end{equation}
where we write $\log$ for the natural logarithm, and where $C>0$ is an
arbitrary constant factor.

\smallskip

In this paper, we attempt to elucidate this surprising connection between the
algebraic structure of the semiring and the information theoretic entropy
by developing a broader theory of thermodynamic semirings.

\section{Preliminary notions}\label{prelimSec}

We introduce here some basic facts that we will need to use in the rest
of the paper.

\smallskip

We start with a warning about notation. Throughout most of the
paper we will work implicitly with $\R^{\min, +} \cup \{ \infty\}$ or 
$\R^{\max, \ast}_{\geqslant 0}$ in mind (note the two are isomorphic under
$- \log$). As such, we will use the notation in one of the two. Which one we use should hopefully be clear from the context. We do this because the first will give expressions looking more like statistical physics equations, and the second will give expressions more similar to the Witt construction in characteristic $p$. We will tend to write $\oplus_{S,T}$ (perhaps with other relevant subscripts) for the Witt addition,
to indicate that it is a modification of the {\em additive} structure of
the semiring, and that it depends on the choice of a binary information
measure (or entropy) $S$ and of a temperature parameter $T$. 
This is motivated by tropical geometry, where
it is customary to denote by $\oplus$ the addition in the tropical
semiring, ie. the minimum, and by $\odot$ the multiplication, the usual addition $+$, see \cite{MaStu}.

\medskip

\subsection{Frobenius in characteristic one}\label{Frob1Sec}

We recall here, from \cite{Co}, \cite{CC}, the behavior of the Frobenius action
in the characteristic one setting.

\smallskip

Let $K$ be a commutative, characteristic one semifield. It is
possible to work in the slightly more general case of multiplicatively
cancellative semirings, but for simplicity we will forsake this generality.
Recall that such a semifield is a set with two associative, commutative binary
operations, $(x, y) \mapsto x + y$ and $(x, y) \mapsto x y$ such that the
second distributes over the first, $0 + x = x$, $0 x = 0$, $1 x = x$, $K$ has
multiplicative inverses, and, importantly, the characteristic one condition
that $1 + 1 = 1$.

\smallskip

The first step in developing an analog of the Witt construction is to
examine the Frobenius map in $K$. 

\begin{lem}\label{oneFrob} (Frobenius)
\begin{equation}\label{Frobn}
 (x + y)^n = x^n + y^n \ \ \  \text{ for every } \ n \in \N. 
\end{equation} 
\end{lem}

\proof
The proof is given in Lemma 4.3 of \cite{Co}, but we recall it here
for the convenience of the readers. One sees from the distributive property that, for
every $m \in \N$, one has $(x + y)^m = \sum_{k = 0}^m x^k y^{m - k}$.
This then gives $(x^n + y^n) (x + y)^{n - 1} = (x + y)^{2 n - 1}$. Since $K$ is
multiplicatively cancellative, this implies \eqref{Frobn}.
\endproof

\medskip

\subsection{Legendre transform}\label{LegendreSec}

As shown in Lemma 4.2 of \cite{Co},
$K$ is endowed with a natural partial ordering $\leqslant$ defined so
that $x \leqslant y \Leftrightarrow x + y = y$. This may seem strange, but
one sees that, over the tropical semifield, $\bT$ this reads $x
\leqslant y \Leftrightarrow \max (x, y) = y$. We give $K$ the order topology
from $\leqslant$. Then multiplication and the Frobenius automorphisms make $K$
a topological $\R_{\geqslant 0}$-module, since the Frobenius is continuous and distributes over the multiplicative structure. When $K =\bT$, this
topology is the standard one on $[0, 1) \cong \R \cup \{- \infty\}$,
with the Frobenius acting by multiplication so that $K$ has the normal vector space
structure.

\smallskip

We say that a function $f : X \rightarrow K$, where $X$ is a convex
subset of a topological $\R_{\geqslant 0}$-module, is convex if, for
every $t \in [0, 1]$, $x_1, x_2 \in X$,
\begin{equation}\label{convexeq}
  f (t x_1 + (1 - t) x_2) \leqslant f (x_1)^t f (x_2)^{1 - t},
\end{equation}
with concavity being defined as convexity of the multiplicative inverse of
$f$. 

Note again that, over $\bT$, this is the normal definition of
convexity. 

\smallskip

We consider also
\begin{equation}\label{episet}
 {\rm epi} f =\{(\alpha, r) \in X \times K | f (\alpha) \leqslant r\}, 
\end{equation} 
called the {\em epigraph} of $f$. This has the following property.

\begin{lem}\label{epilem}
A function $f$ is convex iff the epigraph ${\rm epi} f$ is convex and 
$f$ is closed iff ${\rm epi} f$ is closed.
\end{lem}

\proof The topological $\R_{\geqslant 0}$-module structure on $X
\times K$ is given by the product structure, so the proof follows directly
from the definitions.
\endproof

When $X \subseteq \R_{\geqslant 0}$, we
can define the Legendre transform of $f$ by
\begin{equation}\label{eqLegendre}
  f^{\ast}(x) = \sum_{\alpha \in X} \frac{x^{\alpha}}{f (\alpha)} .
\end{equation}
Note that over $\bT$ this reads
\[ \sup_{\alpha \in X} (\alpha x - f (\alpha)), \]
which is the normal definition of the Legendre transform.

When $X \subseteq K$, we can define the Legendre transform of $f$ by
\begin{equation}
  f^{\ast}(\alpha) = \sum_{x \in X} \frac{x^{\alpha}}{f (x)} .
\end{equation}

\begin{prop}\label{propLegendre}
  The Legendre transform of $f$ is closed and convex.
\end{prop}

\proof
Suppose first $X \subset \R_{\geqslant 0}$. Let $g_{\alpha} (x) = x^{\alpha} / f (\alpha)$, and $g$ be the Legendre
transform of $f$. Then $g$ is the point-wise supremum among the $g_{\alpha}$, so
${\rm epi}\, g = \bigcap_{\alpha \in X} {\rm epi}\, g_{\alpha}$, an intersection
of closed half spaces. Thus, ${\rm epi}\, g$ is closed and convex, so $g$ is
closed and convex, by the Lemma \ref{epilem}. The proof of the opposite case proceeds in precisely the same manner.
\endproof

One then has the following result on Legendre transforms.

\begin{thm}\label{thmLegendre}
  (Fenchel-Moreau) Let $f: X \rightarrow K$, $X \subset \R_{\geqslant 0}$. Then the
  following hold.
\begin{enumerate} 
 \item $f^{\ast \ast}$ is closed and convex and bounded by $f$.
 \item $f^{\ast \ast} = f$ iff $f$ is closed and convex.
\end{enumerate} 
\end{thm}

\proof The function $f^{\ast \ast}$ is convex and closed by Lemma \ref{epilem}. 
We also see that
\[ x^{\alpha} / f^{\ast} (x) \leqslant x^{\alpha} / (x^{\alpha} / f (\alpha)) =
f (\alpha), \] 
so taking a $\sup_{x \in X}$ of both sides yields $f^{\ast \ast}
\leqslant f$. To prove the second fact, it suffices to show that, if $f$ is closed,
convex, and finite, then $f \leqslant f^{\ast \ast}$. Define the
subdifferential $\partial f (\alpha)$ of $f$ at $\alpha$ by
\[ \partial f (\alpha) =\{x \in R | f (\beta) \geqslant f (\alpha) x^{\beta -
   \alpha} \forall \beta \in X\}. \]
We consider the set-valued map $\alpha \mapsto \partial f (\alpha)$. To invert
this map is to find $\alpha (x) = \alpha$ such that $x \in \partial f
(\alpha)$. We see $f^{\ast} (x) = x^{\alpha (x)} / f (\alpha (x))$. Thus, the
subdifferential is the proper analog in this case for the derivative. When $f$
is closed and convex, $\partial f (\alpha)$ is nonempty, so let $x \in
\partial f (\alpha)$. Then we have
\[ 1 / f^{\ast} (x) \geqslant f (\alpha) / x^{\alpha} \Rightarrow f (\alpha)
   \leqslant x^{\alpha} / f^{\ast} (x)) \leqslant f^{\ast \ast} (\alpha) \]
for every $\alpha$, proving the theorem.
\endproof

\smallskip

This is a simple translation of the well-known Legendre transform machinery into characteristic one semirings. The idea is that since we can define a real topological vector space structure on $K$ using the multiplication as addition and the Frobenius map as scalar multiplication (with negative reals having a well-defined action since $K$ has multiplicative inverses), we have enough structure to do convex analysis. The point is that for concave or convex $f$, the above sums are invertible in $K$. From now on, any semifield satisfying the assumptions necessary for this section will be called ``suitably nice".

\subsection{Witt ring construction in characteristic one}

We recall here the main properties of the characteristic one 
analog of the Witt construction \cite{Co}, \cite{CC}, which is the 
starting point for our work. We formulate
it here in terms of a general information measure $S$, whose properties 
we will find are related to the algebraic properties of the semiring.

\smallskip

Let $w : [0, 1] \rightarrow K$ be continuous under the order
topology, and consider, for each $x, y \in K$,
\begin{equation}\label{wsum}
 x \oplus_w y = \sum_{\alpha \in I} w (\alpha) x^{\alpha} y^{1 - \alpha} . 
\end{equation}

\smallskip

Connes and Consani considered the above expression for continuous $w
(\alpha) \geq 1$ and found in \cite{Co}, \cite{CC} that $\oplus_w$ is commutative, associative and has
identity $0$ if and only if $w (\alpha) = \rho^{{\rm Sh} (\alpha)}$ for some $\rho \in
K$ greater than one. 

\smallskip

For simplicity and clarity of intention, we will write
$\rho = e^T$ for some $T \geqslant 0$  to suggest $T$ behaves like a temperature parameter. In all the arguments that follow, one could
replace $e^T$ by $\rho$ again and be fine over the more general semifields.

\smallskip

Correspondingly, we are going to restrict our attention to sums of the form
\begin{equation}\label{sumTS}
  x \oplus_S y := \sum_{\alpha \in I} e^{T S (\alpha)} x^{\alpha} y^{1 - \alpha}
\end{equation}
where $S$ will be interpreted as an entropy function. In particular, we
assume $S$ is concave and closed, so that $e^{- T S (\alpha)}$ is convex and
closed, and we can use the Legendre transform machinery developed 
in \S \ref{LegendreSec}.

\smallskip

We can then formulate the result of \cite{CC} on the characteristic one Witt
construction in the following way.

\begin{thm}\label{CConeWitt}
Suppose $S : I \rightarrow \R_{\geqslant 0}$ is concave and closed.
  The following hold.
  
 \begin{enumerate} 
  \item $x \oplus_S y = y \oplus_S x$ $\forall x, y \in K$ iff $S (\alpha) = S (1 -
  \alpha)$.  
  \item $0 \oplus_S x = x$ $\forall x \in K$ iff $S (0) = 0$.  
  \item $x \oplus_S 0 = x$ $\forall x \in K$ iff $S (1) = 0$.  
  \item $x \oplus_S (y \oplus_S z) = (x \oplus_S y) \oplus_S z$ $\forall x, y, z \in K$
  iff $S (\alpha \beta) + (1 - \alpha \beta) S ( \frac{\alpha (1 - \beta)}{1 -
  \alpha \beta}) = S (\alpha) + \alpha S (\beta)$.
\end{enumerate}
\end{thm} 
 
 \proof 
The argument is given in \cite{Co} in a more general form 
 applicable to a binary operation as in \eqref{charonesums}, but we give the explicit proof here to show the machinery.
 
 \smallskip
 
(1) Since $S$ is concave and closed, $e^{- T S}$ is convex and closed
(in the generalized sense of \eqref{convexeq}), as is $z^{L (\alpha)}$ for any linear
function $L (\alpha)$ and $z \in K$. We also see that products of convex and
closed functions are convex and closed, so $y^{\alpha - 1} e^{- T S (\alpha)}$
and $y^{\alpha - 1} e^{- T S (1 - \alpha)}$ are each convex and closed. We see that
$x \oplus_S y = y \oplus_S x$ iff
\[ \sum_{\alpha \in I} \frac{x^{\alpha}}{y^{\alpha - 1} e^{- T S (\alpha)}} =
   \sum_{\alpha \in I} \frac{x^{\alpha}}{y^{\alpha - 1} e^{- T S (1 -
   \alpha)}} . \]
We recognize the Legendre transform of closed convex functions, which is invertible by the Fenchel-Moreau theorem above. Thus, the summands must be
equal, so $S (\alpha) = S (1 - \alpha)$. The converse is obvious.

\smallskip

(2) First note that, when $\alpha \neq 0$, for every $x$, $0^{\alpha} x^{1 -
\alpha} = 0$ and $e^{T S(0)}\geq 0$, so the supremum occurs at $\alpha = 0$. Therefore, we have $0
\oplus_S x = e^{T S (0)} x$.

\smallskip

(3) Similarly, this supremum occurs at $\alpha = 1$, so $x \oplus_S 0 = e^{T S
(1)} x$.

\smallskip

(4) As in fact {\em 1}, we see $x \oplus_S (y \oplus_S z) = (x \oplus_S y) \oplus_S z$ iff
\[ \sum_{\alpha, \beta \in I} \frac{x^{\alpha \beta}}{y^{\alpha (\beta - 1)}
   z^{\alpha - 1} e^{- T (S (\alpha) + \alpha S (\beta))}} = \sum_{u, v \in I}
   \frac{x^u}{y^{v (u - 1)} z^{(v - 1) (1 - u)} e^{- T (S (u) + (1 - u) S
   (v))}} . \]
Identifying powers and inverting the Legendre transform yields the condition.
The converse is immediate. 
\endproof

\smallskip

We hold off discussing the fact that the Shannon entropy ${\rm Sh}$ is the only function $S$
satisfying all of these properties until \S
\ref{Shaxiomsec} below, where we develop the information theoretic interpretation of these axioms.

\section{Axioms for Entropy Functions}\label{Shaxiomsec}

It is well known that the Shannon entropy admits an axiomatic characterization
in terms of the Khinchin axioms \cite{Khin}. These are usually stated as follows for
an information measure  $S(p_1,\ldots, p_n)$:
\begin{enumerate}
\item (Continuity) For any $n \in \N$, the function $S(p_1,\ldots, p_n)$ is 
continuous with respect to $(p_1,\ldots, p_n)$ in the simplex 
$\Delta_n=\{ p_i\in \R_+, \, \sum_i p_i =1\}$;
\item (Maximality) Given $n\in \N$ and $(p_1,\ldots, p_n)\in \Delta_n$, the
function $S(p_1,\ldots, p_n)$ has its maximum at the uniform distribution $p_i=1/n$
for all $i=1,\ldots,n$,
\begin{equation}\label{Khin2}
  S(p_1,\ldots, p_n) \leq S(\frac{1}{n}, \ldots, \frac{1}{n}), \ \ \  \forall (p_1,\ldots, p_n)\in \Delta_n; 
\end{equation}  
\item (Additivity) If $p_i=\sum_{j=1}^{m_i} p_{ij}$ with $p_{ij}\geq 0$, then
\begin{equation}\label{Khin3}
 S(p_{11},\ldots, p_{nm_n}) =
 S(p_1,\ldots, p_n) + \sum_{i=1}^n p_i S(\frac{p_{i1}}{p_i}, \ldots, \frac{p_{im_i}}{p_i});
\end{equation}
\item (Expandability) Embedding a simplex $\Delta_n$ as a face inside a simplex $\Delta_{n+1}$
has no effect on the entropy,
\begin{equation}\label{Khin4}
S(p_1,\ldots,p_n,0)= S(p_1,\ldots, p_n).
\end{equation} 
\end{enumerate}
It is shown in \cite{Khin} that there is a unique 
information measure $S(p_1,\ldots,p_n)$ (up to a multiplicative constant $C>0$)  
that satisfies these axioms and it is
given by the Shannon entropy 
\begin{equation}\label{ShKhin}
S(p_1,\ldots,p_n)={\rm Sh}(p_1,\ldots,p_n):= - C \sum_{i=1}^n p_i \log p_i ,
\end{equation}

\medskip

We focus now on the $n=2$ case, which means that 
we are only looking at $S(p):=S(p,1-p)$ instead of the more
general $S(p_1,\ldots,p_n)$. In other words, we are only considering the information theory of binary random variables.
In this case, we describe here an axiomatic formulation 
for the Shannon entropy based
on properties of binary ``decision machines".
We return to discuss the
more general $n$-ary case in  \S \ref{operadsec} below.

\smallskip

A decision machine is a measurement tool which may only distinguish between two possible states of
a discrete random variable; machines that can only answer ``yes'' or ``no''.
We would like to measure the average change in uncertainty after a
measurement, which is how we define the entropy associated with a random
variable. Let $X$ be a binary random variable, $S (X)$ the change in entropy
after measuring $X$. All information is created equal, so $S (X)$
should only depend on the probability of measuring a certain value of $X$ and should do so continuously.

\smallskip

\begin{enumerate}
\item   (Left Identity) $S (0) = 0$. 
\item   (Right Identity) $S (1) = 0$.
\item   (Commutativity) $S (p) = S (1 - p)$.
\item   (Associativity) $S (p_1) + (1 - p_1) S ( \frac{p_2}{1 - p_1}) = S (p_1 +
  p_2) + (p_1 + p_2) S ( \frac{p_1}{p_1 + p_2})$.
\end{enumerate}

\smallskip

The identity axioms claim that trivial measurements 
give trivial information.

\smallskip

The commutativity axiom claims that questions have the same information as their negative.

\smallskip

The associativity axiom claims a certain equivalence of guessing strategies, which will be a key observation in our explanation of the characteristic one Witt construction. 
If instead of a binary random variable, we want to measure a ternary random variable $X$ which may take
values $X \in \{x_1, x_2, x_3 \}$ with corresponding probabilities $p_1, p_2,
p_3$, we can still determine $X$ by asking yes-or-no questions. We can first
ask ``is $X = x_1$?'' If the answer is no (which occurs with probability $p_2
+ p_3$), we then ask ``is $X = x_2$?'' This corresponds to an average change
in uncertainty $S (p_1) + (p_2 + p_3) S ( \frac{p_2}{p_2 + p_3})$. However, we
could have asked ``is $X = x_1$ or $x_2$ ?'' followed by ``is $X =
x_1$?'' and in the end received the same data about $X$. Associativity asserts these two
should be equal, hence we have the axiom as stated above. 

\smallskip

The names of the axioms in the above list are chosen to suggest the corresponding
algebraic properties, as we see in Theorem \ref{algSthm} below. 
In fact, we find that these algebraicly motivated axioms are equivalent to the Khinchin axioms. 

\begin{thm}\label{thmSaxioms}
There is a unique function (up to a multiplicative constant $C>0$) 
satisfying all of the axioms above, namely the Shannon entropy
\begin{equation}\label{Shp}
  {\rm Sh} (p) = - C (p \log p + (1 - p) \log (1 - p)) .
\end{equation}  
\end{thm}

\proof The result follows either by checking directly the equivalence of the
commutativity, identity and associativity axioms with the Khinchin axioms,
or else by proceeding as in Theorem 5.3 of \cite{Co}. We prove it here
by showing that one obtains the Khinchin axioms for entropy.

\smallskip

Suppose $S$ satisfies all the conditions above. Define $S_n : \Delta_{n
- 1} \rightarrow \R_{\geqslant 0}$ by
\begin{equation}\label{SnSdef}
 S_n (p_1, \ldots, p_n) = \sum_{1 \leqslant j \leqslant n - 1} (1 - \sum_{1
   \leqslant i < j} p_i) S ( \frac{p_j}{1 - \sum_{1 \leqslant i < j} p_i}) .
\end{equation}

\begin{lem}
  $S_n$ is symmetric.
\end{lem}

\proof
Suppose we interchange the terms $p_k$ and $p_{k + 1}$, where $k < n - 1$.
This only affects the terms $k$th and $(k + 1)$th terms, so we must show
\[ T = (1 - \sum_{i < k} p_i) S ( \frac{p_k}{1 - \sum_{i < k} p_i}) + (1 -
   \sum_{i < k + 1} p_i) S ( \frac{p_{k + 1}}{1 - \sum_{i < k + 1} p_i}) \]
is symmetric. Write $\beta = 1 - \sum_{i < k} p_i$, $a = p_k / \beta$, $b =
p_{k + 1} / \beta$. We see $\beta$ is invariant under this permutation, and
\[ T = \beta \left( S (a) + (1 - a) S (b / (1 - a))\right) . \]
Permuting $p_k$ and $p_{k + 1}$ interchanges $a$ and $b$, and so $T$ is
invariant by the associativity condition. Interchanging $p_{n - 1}$ and $p_n$
only affects the last term, and it is easy to see it affects it like $S
(\alpha) \mapsto S (1 - \alpha)$, so invariance follows from commutativity.
These transpositions generate the symmetric group ${\rm Sym}_n$, so $S_n$ is
symmetric.
\endproof

From this lemma and the definition we see the following holds.

\begin{lem}
  Let $(J_k)_{1 \leqslant k \leqslant m}$ be a partition of $\{p_1, \ldots,
  p_n \}$ and let $S_n$ be defined as in \eqref{SnSdef}. Then we have
  \[ S_n (p_1, \ldots, p_n) = S_m (q_1, \ldots, q_m) + \sum_{1 \leqslant k
     \leqslant m} S_{|J_k |} (J_k / q_k), \]
  where $q_k = \sum_{p_{} \in J_k} p$, so $J_k / q_k$ is a $|J_k |$-ary
  probability distribution.
\end{lem}

These lemmas take care of the third Khinchin axiom, and with the identity property also take care of the fourth. We assumed at the outset $S$ was continuous, so it follows from the definition $S_n$ is continuous, which is the first axiom. What remains is the second axiom, which we
write here in terms of information (concave) rather than entropy (convex).

\smallskip

\begin{lem}
  $S_n$ is concave for all $n$.
\end{lem}

\proof
We  proceed by induction on $n$. We have already assumed $S_2 = S$ is concave, so suppose $S_n$ is concave for some $n \geq 2$. Note that for continuous $f$, concavity follows from $f(\frac{x+y}{2})\geq \frac{f(x)+f(y)}{2}$. Thus we consider, for some $(p_i), (q_i) \in \Delta_{n+1}$
\[S_{n+1}(\frac{p_1}{2}+\frac{q_1}{2}, ..., \frac{p_{n+1}}{2}+\frac{q_{n+1}}{2}).\]
By the previous lemma, this equals
\[S_{n}(\frac{p_1+p_2}{2}+\frac{q_1+q_2}{2}, \frac{p_3}{2}+\frac{q_3}{2}, ... , \frac{p_{n+1}}{2}+\frac{q_{n+1}}{2}) \] \[+\frac{p_1+p_2+q_1+q_2}{2}S(\frac{p_1+q_1}{p_1+p_2+q_1+q_2}).\]
By the inductive hypothesis we then have
\begin{align*}
S_{n+1}(...) \geq & \frac{1}{2}S(p_1+p_2, ...,p_{n+1})+\frac{1}{2}S(q_1+q_2,...,q_{n+1}) \\
+& \frac{p_1+p_2+q_1+q_2}{2}S(\frac{p_1+q_1}{p_1+p_2+q_1+q_2}).\end{align*}
We see that $$\frac{p_1+q_1}{p_1+p_2+q_1+q_2}=\frac{p_1+p_2}{p_1+p_2+q_1+q_2}\frac{p_1}{p_1+p_2}+\frac{q_1+q_2}{p_1+p_2+q_1+q_2}\frac{q_1}{q_1+q_2}$$ and $$\frac{p_1+p_2}{p_1+p_2+q_1+q_2}+\frac{q_1+q_2}{p_1+p_2+q_1+q_2}=1,$$ 
so by the concavity of $S$ we have
\begin{align*}
S_{n+1}(...) \geq & \frac{1}{2}S(p_1+p_2, ...,p_{n+1})+\frac{1}{2}S(q_1+q_2,...,q_{n+1}) \\ +& \frac{p_1+p_2}{2}S(\frac{p_1}{p_1+p_2})+\frac{q_1+q_2}{2}S(\frac{q_1}{q_1+q_2}), \end{align*}
from which concavity of $S_{n+1}$ follows by the previous lemma.
\endproof

Since $S_n$ is concave, it has a unique maximum, and since it is symmetric, this maximum occurs at $S_n(\frac{1}{n},...,\frac{1}{n})$, implying the second Khinchin axiom. This then completes the proof of Theorem \ref{thmSaxioms}
\endproof

A reformulation of the Khinchin axioms for Shannon entropy
more similar to the commutativity, identity and associativity axioms
considered here was described in Faddeev's \cite{Fad}. For different
reformulations of the Khinchin axioms see also \cite{Csi}.

\section{Thermodynamic semirings}\label{thermoSec}

We now consider more general thermodynamic semirings. The following
definition describes the basic structure.

\begin{defn}\label{thermosemidef}
  A thermodynamic semiring structure over $K$, written like $\R^{\min, +} \cup \{\infty\}$, 
  is a collection of binary operations $\oplus_{S,T} : K \times
  K \rightarrow K$ indexed by $T \in \R \cup \{\infty\}$ and defined by an
  information measure, $S : [0, 1] \rightarrow \R$ according to 
\begin{equation}\label{plusTgen}  
  x \oplus_{S,T} y = \min_{p \in [0, 1] \cap \Q} (p x + (1 - p) y - T S
  (p)).
\end{equation}  
\end{defn}

It is often convenient to consider the elements of the semiring as functions of $T$, with the
operation $\oplus_S$ defined pointwise by $\oplus_{S,T}$. We call this ring $R$, inspired the 
$p$-typical Witt notation. Indeed in \cite{CC}, \cite{Co}, $R$ is seen as the Witt ring over $K$, 
with evaluation at $T=0$ over giving the residue morphism $R \rightarrow K$. We then see that the Teichm\"uller lifts 
should be the constant functions, and $T$ should play the role of the exponent of $p^n$ in considering field extensions.

We then have the following general properties, as in Theorem \ref{thmSaxioms}
above (Theorem 5.2 of \cite{Co}):

\begin{thm}\label{algSthm}
Let $x\oplus_{S,T} y$ be a thermodynamic semiring structure on a suitably nice characteristic one semifield, $K$, defined as in \eqref{plusTgen}.
Then the following holds.
\begin{enumerate}  
 \item  $x \oplus_{S,T} y = y \oplus_{S,T} x$ iff $S$ is commutative.
 \item  $0 \oplus_{S,T} x = x$ iff $S$ has the left identity property.
 \item  $x \oplus_{S,T} 0 = x$ iff $S$ has the right identity property.
 \item  $x \oplus_{S,T} (y \oplus_{S,T} z) = (x \oplus_{S,T} y) \oplus_{S,T} z$ iff $S$ is associative.
\end{enumerate} 
\end{thm}

\proof The case of commutativity and of the identity axioms are 
obvious. For associativity we have
\[ x \oplus_{S,T} (y \oplus_{S,T} z) = x \oplus_{S,T} \min_p (p y + (1 - p) z - T S (p)) \]
\[ = \min_q (q x + (1 - q) \min_p (p y + (1 - p) z - T S (p)) - T S (q)) \]
\[ = \min_{p, q} (q x + p (1 - q) y + (1 - q) (1 - p) z - T (S (q) + (1 - q) S
(p))) \]
\[ = \min_{p_1 + p_2 + p_3 = 1} (p_1 x + p_2 y + p_3 z - T (S (p_1) + (1 - p_1) S
( \frac{p_2}{1 - p_1}))) \]
while
\[  (x \oplus_{S,T} y) \oplus_{S,T} z = \min_p (p x + (1 - p) y - T S (p)) \oplus_S z \]
\[ = \min_{p, q} (p q x + q (1 - p) y + (1 - q) z - T (q S (p) + S (q)) \]
\[ = \min_{p_1 + p_2 + p_3 = 1} (p_1 x + p_2 y + p_3 z - T (S (p_1 + p_2) + (p_1
+ p_2) S ( \frac{p_1}{p_1 + p_2})). \]
We see that the two ways of summing three quantities corresponds to the two
ways of measuring a ternary random variable with decision machines. The
equivalence is now obvious.
\endproof

Most information measures are commutative, though a non-commutative example in \S \ref{KLsec} below. We discuss
in \S \ref{statmechSec} some physical reasons why commutativity
is more automatic in this context than associativity.

\medskip

One then sees by direct inspection that, in the case of the Shannon entropy
one has the following form of the thermodynamic semiring structure.

\begin{prop}\label{propShsum}
  When $S$ is the Shannon entropy, ${\rm Sh}$, then
  \begin{equation}\label{plusShlog}
  x \oplus_{{\rm Sh},T} y = - T \log (e^{- x / T} + e^{- y / T})
  \end{equation}
   over  $\R^{\min, +} \cup \{\infty\}$, while over $\R^{\max,  \ast}_{\geqslant 0}$ it is
  \begin{equation}\label{plusShT}
  x \oplus_{{\rm Sh},T} y = (x^{1 / T} + y^{1 / T})^T  .
  \end{equation}
\end{prop}

\smallskip

Notice that the semiring $\R^{\max,  \ast}_{\geqslant 0}$ is isomorphic to
the semiring $\R^{\min, +} \cup \{\infty\}$, under the $-\log$ mapping, so that
\eqref{plusShT} is simply obtained from \eqref{plusShlog} in this way.

\medskip

In this case, the parameter $T$ corresponds to the parameter $h$ of
Maslov dequantization (see the comments in \S \ref{tropicalsec}). The
semifields obtained in this way are known as the Gibbs--Maslov semirings
and the subtropical algebra (see \cite{LoMaQua}, \cite{Litv2}).

\medskip

One can extend the notion of thermodynamic semiring to include a class
of semirings of functions which we will be considering in the following.
Just as in the case of a ring $R$ and a parameter space $X$, one can 
endow the set of functions from $\Xi$ to $R$ with a ring structure, 
by pointwise operations, one can proceed similarly with a semiring.
Moreover, in the case of a thermodynamic semiring structure, it is
especially interesting to consider cases where the pointwise operation
$\oplus_{S,T}$ depends on the point in the parameter space through a
varying entropy function $S=S_\eta$, for $\eta\in \Xi$.

\begin{defn}\label{funthermosemi}
Let $\Xi$ be a compact Hausdorff space and let $S=(S_\eta)$ be a 
family of information measures depending continuously on the
parameter $\eta\in \Xi$. Let $K=\R^{\min, +} \cup \{\infty\}$. A
thermodynamic semiring structure on the space of functions 
$C(X,R)$ is given by the family of pointwise operations
\[ x(\eta) \oplus_{S_\eta,T} y(\eta) =\min_{p\in [0,1]\cap \Q} (px(\eta) +(1-p)y(\eta) -T S_\eta(p)). \]
\end{defn}

The properties of Theorem \ref{algSthm} extend to this case. 
We will return to this more general setting in \S \ref{KLsec} below.

\medskip

As we discuss in the following sections, more general entropy
functions (which include the special cases of R\'enyi entropy, 
Tsallis entropy and Kullback--Leibler divergence, as well as 
the more general categorical and operadic setting developed 
in \S \ref{operadsec}) give rise to thermodynamic algebraic
structures that are neither commutative nor associative. 
We will continue to use the terminology ``semiring", although
(as the referee pointed out to us) the term ``algebra",
in the sense of the theory of universal algebra, would be
more appropriate.

\section{Statistical mechanics}\label{statmechSec}

Before we move on to see explicit examples of thermodynamic
semirings besides the original one based on the Shannon entropy
considered already in \cite{Co} and \cite{CC}, we give in this section a 
physical interpretation of the algebraic
structure of thermodynamic semirings in terms of statistical mechanics. 
This interpretation is a generalization of thermodynamic interpretations 
of max-plus formulas found in \cite{Quad}.

\medskip

When $K =\R^{\max, \ast}_{\geqslant 0}$, we can write the thermodynamic semiring
operations in the form
\[ x \oplus_{\rho,S} y = \max_p (\rho^{S (p)} x^p y^{1 - p}). \]
In particular, when we set $\rho = e^{k_B T}$, this reads
\[ \max_p (e^{k_B T S(p)+ p \log x + (1 - p) \log y}) .\]

We recognize this as $e^{- F_{\rm eq}} = Z$, where $F_{\rm eq}$ is the
equilibrium value of the free energy of a system at temperature $T$,
containing a gas of particles with chemical potentials $\log x$ and $\log y$,
and Hamiltonian
\begin{equation}\label{Hamiltonian}
 \mathcal{H}= p \log x + (1 - p) \log y,  
\end{equation} 
where $p$ is now thought of as a mole fraction, and $Z$ is its partition
function.

\smallskip

Indeed, the semirings $\R^{\max, \ast}_{\geqslant 0}$ and
$\R^{\min, +} \cup \{ \infty\}$ are isomorphic by $- \log$, and this
gives
\[ \log x \oplus_{S, k_B T} \log y = \min_p (p \log x + (1 - p) \log y - k_B T S (p)), \]
which is the equilibrium free energy described above. We note also that the calculated form of the thermodynamic semiring for Shannon entropy, that is
\[x\oplus_{Sh}y = - T \log(e^{-x/T}+e^{-y/T}).\]
In it, we recognize precisely the partition sum of a two state system with energies $x$ and $y$. We thus consider members of the Witt ring $R$ (see \S \ref{thermoSec}) to be temperature dependent chemical potentials.

\medskip

In a gas system with a single type of particle, the free energy is precisely the
chemical potential. The mixing of these gases gives a new free energy dependent on the entropy function. We then replace this mixture with a ``particle" whose chemical potential is the equilibrium free energy per particle of the previous mixture. This gives a monoid structure on the space of chemical potentials. When we consider mixing in arbitrary thermodynamics, ie. with non-Boltzmann counting, we
have the possibility of mixing to be {\em non-associative}. With this
interpretation, however, we would not expect the mixing process to ever be
{\em non-commutative}, so the lack of associativity has a more direct and natural physical interpretation than the lack of commutativity for
thermodynamic semirings. We imagine multiplication to be a sort of bonding of gases, where chemical potentials add together.
\medskip

We see the dynamics of this mixing process is determined, both physically and algebraicly, by the entropy
function and the ambient temperature. At zero temperature, the mixture is always entirely
composed of the particle with the least chemical potential. This corresponds
to $\R^{\min, +}$ and indeed evaluation at zero temperature gives us the residue morphism $R \rightarrow K$. When the entropy function is the Shannon entropy,
we get the normal thermodynamical mixing, see \S 8.5 of \cite{Denb}. We can say, therefore, that the Witt construction is, in a sense, giving thermodynamics to this system. Note that in \eqref{charonesums} this construction is seen giving an inverse to Maslov dequantization, pointing out an interesting link between quantum mechanics and thermodynamics.

\medskip

The mixing entropy for chemical systems based on the Boltzmann--Gibbs
statistical mechanics and the Shannon entropy
function (as in \S 8.5 of \cite{Denb} for instance) works well to describe 
systems that are {\em ergodic}.
If a system is nonergodic (that is, time averages and phase space averages differ),
then the counting involved in bringing two initially separated systems into
contact will not follow the normal Boltzmann rules. As a result, Shannon
entropy will not behave extensively in these systems. This typically occurs in
physical systems with strong, long-range coupling and in systems with metastable states or exhibiting power law
behavior. In such systems, maximizing the Shannon entropy functional (subject to the
dynamical constraints of the system) does not produce the correct
metaequilibrium distribution, see for instance \cite{Tsallis} and other essays
in the collection \cite{GeTsa}.

\smallskip

A broad field of nonextensive statistical mechanics for such
systems has been developed (see \cite{Tsallis} for a brief introduction),
where, under suitable conditions, one can calculate a ``correct'' entropy functional
corresponding to the system at hand. These entropy functionals are typically
characterized by some axiomatic properties that describe their behavior.
For instance, if we have two initially
independent systems $A, B$ and bring them together to form a combined
system denoted $A \star B$, one may require that 
$S (A \star B) = S (A) + S (B)$ (extensive).  This leads to forms of entropy such as
the R\'enyi entropy \cite{Renyi}, generalizing the original Shannon case, while maintaining
the extensivity over independent systems. One may also have explicit $q$-deformations of the extensivity condition, for example $S_q(A\star B)= S_q(A)+ S_q(B)+ (1-q) S_q(A) S_q(B)$ for independent
systems. This leads to forms of entropy such as the Tsallis entropy \cite{Tsallis2}.

\smallskip

When we consider different kinds of entropy functions in this way, we can
look at the algebraic properties of the corresponding thermodynamic semirings.
These will encode the information about the amount of nonextensivity and
nonergodicity of the system giving rise to the corresponding entropy function $S$. We can imagine non-associativity of mixing as a toy model of meta-equilibrium states where we known the entropy beforehand.
We can also use thermodynamic semirings to encode relative entropies and
analyze its behavior over a space of parameters through the algebraic 
properties of the semiring.

\medskip

Relations between idempotent semifields and statistical mechanics were
also considered in \cite{ItMik}, \cite{Kapr}, \cite{Quad}.

\section{The R\'enyi entropy}\label{RenyiSec}

We now look at other important examples of entropy functions and we investigate
how the corresponding algebraic properties of the associated thermodynamic
semiring detect the properties of the entropy function as an information measure.

\smallskip

A first well known case of an entropy function which is
a natural generalization of the Shannon entropy: the R\'enyi entropy, \cite{Renyi}. 
This is a one-parameter family ${\rm Ry}_\alpha$ of information
measures defined as 
\begin{equation}\label{Renyientropy}
{\rm Ry}_\alpha (p_1,\ldots,p_n) := \frac{1}{1-\alpha} \log \left( \sum_i p_i^\alpha \right),
\end{equation} 
so that the limit
\begin{equation}\label{Renyilim}
\lim_{\alpha \to 1} {\rm Ry}_\alpha (p_1,\ldots,p_n) = {\rm Sh}  (p_1,\ldots,p_n) 
\end{equation}
recovers the Shannon entropy. The R\'enyi entropy has a broad range of
applications, especially in the analysis of multifractal systems \cite{BeSch},
while a statistical mechanics based on the R\'enyi entropy is described
in \cite{LeMeSi}. 

\smallskip

The R\'enyi entropy also has an axiomatic characterization, where
one weakens the Khinchin additivity axioms to a form that only requires
additivity of the information entropy for independent subsystems, while
keeping the other three axioms unchanged, \cite{Renyi2}.
For our version of the axioms, formulated in terms of decision machines, 
this means that the associativity axiom no longer holds. 

\smallskip 

\begin{lem}\label{Renyiassoc}
The lack of associativity of $x\oplus_S y$, when $S={\rm Ry}_\alpha$ is
the Renyi entropy
\begin{equation}\label{Renyip}
{\rm Ry}_\alpha(p)= \frac{1}{1-\alpha} \log (p^\alpha +(1-p)^\alpha),
\end{equation}
is measured by the transformation $(p_1,p_2,p_3)\mapsto (p_3,p_2,p_1)$.
\end{lem}

\proof We have
\[ {\rm Ry}_\alpha(p_1) + (1-p_1) {\rm Ry}_\alpha (\frac{p_2}{1-p_1}) = \]
\[ \frac{1}{1-\alpha} \left( \log (p_1^\alpha +(1-p_1)^\alpha) + (1-p_1) \log \left(
(\frac{p_2}{1-p_1})^\alpha + ( \frac{ 1-p_1-p_2 }{1-p_1})^\alpha \right)\right) = \]
\[  \frac{1}{1-\alpha}  \log \left( (p_1^\alpha +(1-p_1)^\alpha) \frac{  
 (\frac{p_2}{1-p_1})^\alpha + ( \frac{ p_3 }{1-p_1})^\alpha }{  ( (\frac{p_2}{1-p_1})^\alpha + ( \frac{ p_3 }{1-p_1})^\alpha )^{p_1}   } \right) = \]
\[  \frac{1}{1-\alpha}  \log \left( ( \frac{ p_1 p_2 }{1-p_1} )^\alpha +  ( \frac{ p_1 p_3 }{1-p_1} )^\alpha
+ p_2^\alpha + p_3^\alpha \right) - \frac{p_1}{1-\alpha} \log \left( (\frac{p_2}{1-p_1})^\alpha + ( \frac{ p_3 }{1-p_1})^\alpha \right) \]
\[ \frac{1}{1-\alpha} \log \left( \frac{ (p_2^\alpha + p_3^\alpha) (p_1^\alpha +(1-p_1)^\alpha) }
{(1-p_1)^\alpha} \right) -\frac{p_1}{1-\alpha} \log \left( \frac{(p_2^\alpha + p_3^\alpha)}{(1-p_1)^\alpha} \right) \]
\[  = \frac{1}{1-\alpha} \left( (1-p_1) \log( p_2^\alpha + p_3^\alpha ) + \log  (p_1^\alpha +(1-p_1)^\alpha) - \alpha (1-p_1) \log(1-p_1)\right). \]
On the other hand, we have
\[ {\rm Ry}_\alpha(p_1+p_2) +(p_1+p_2) {\rm Ry}_\alpha(\frac{p_1}{p_1+p_2}) = \]
\[  {\rm Ry}_\alpha(1-p_3) +(1-p_3) {\rm Ry}_\alpha(\frac{p_1}{1-p_3}) =
{\rm Ry}_\alpha(p_3) + (1-p_3) {\rm Ry}_\alpha(\frac{p_1}{1-p_3}) = \]
\[  \frac{1}{1 - \alpha} \log (( \frac{p_1 p_3}{1 - p_3})^{\alpha} + ( \frac{p_2
p_3}{1 - p_3})^{\alpha} + p_1^{\alpha} + p_2^{\alpha}) - \frac{p_3}{1 -
\alpha} \log (( \frac{p_1}{1 - p_3})^{\alpha} + ( \frac{p_2}{1 -
p_3})^{\alpha}) \]
\[ = \frac{1}{1-\alpha} \left( (1-p_3) \log( p_2^\alpha + p_1^\alpha ) + \log  (p_3^\alpha +(1-p_3)^\alpha) - \alpha (1-p_3) \log(1-p_3)\right). \]
So the failure of associativity is corrected by mapping $(p_1,p_2,p_3)\mapsto (p_3,p_2,p_1)$. In fact, this holds for any commutative $S$.

In a commutative non-associative semiring $K$, the lack of associativity is corrected by the
morphism
$$\xymatrix{ K\otimes K \otimes K    \ar@{->}[rr]^{A} \ar@{->}[d]_{\oplus_w\otimes 1} & &
K\otimes K \otimes K \ar@{->}[d]^{1\otimes \oplus_w}  \\ 
K\otimes K \ar@{->}[r]^{\oplus_w} & K & K\otimes K \ar@{->}[l]_{\oplus_w} } $$
which makes the diagram commutative, and which 
is simply given by $A(x\otimes y\otimes z)=z\otimes y \otimes x$.
This is exactly the transformation $(p_1,p_2,p_3)\mapsto (p_3,p_2,p_1)$,
as these correspond to $p_1=sr$, $p_2=s(1-r)$ and $p_3=1-(p_1+p_2)$ in the
associativity constraints. Thus, the transformation $(p_1,p_2,p_3)\mapsto (p_3,p_2,p_1)$
is exactly the one that identifies $w(s) w(r)^s$ with $w(sr) w(s(1-r)/1-sr)^{1-sr}$.
\endproof 

We will show in \S \ref{assocSec} below that one can introduce a more refined
notion of {\em successor function} for thermodynamic semirings, which encodes useful
information on the algebraic structure of the semiring, including the lack of
associativity, and on the thermodynamical properties of the entropy function.

\section{The Tsallis entropy}\label{tsallisSec}

The Tsallis entropy \cite{Tsallis2} is a well-studied generalization of
Shannon entropy, currently finding application in the statistical mechanics of
nonergodic systems, \cite{GeTsa}.  It is defined by  
\begin{equation}\label{Tsallis}
{\rm Ts}_{\alpha} (p) = \frac{1}{\alpha - 1} (1 - p^{\alpha} - (1 -p)^{\alpha}).
\end{equation}
(A slightly more general form will be analyzed in \S \ref{WittTsallisSec} below,
see \eqref{phiTsallis}.)

\medskip

The basic characterizing feature of the Tsallis entropy is the fact that
the extensive property (additivity on independent subsystems) typical
of the Shannon and R\'enyi entropies is replaced by a nonextensive
behavior. This corresponds, algebraically, to replacing an exponential
function (or a logarithm) with an $\alpha$-deformed exponential (or logarithm),
see \S 2.1 of \cite{Tsallis}, so that the usual Boltzmann principle $S= k \log W$
of statistical mechanics is replaced by its deformed version $S_\alpha= k \log_\alpha W$,
where $\log_\alpha(x)=(x^{1-\alpha}-1)(1-\alpha)$. Thus, instead of additivity
$S(A\star B)=S(A)+S(B)$ on the combination of independence systems, one
obtains $S_\alpha(A\star B)=S_\alpha(A)+S_\alpha(B)+(1-\alpha) S_\alpha(A)S_\alpha(B)$.
An axiomatic characterization of the Tsallis entropy is described in
 \cite{Furu}, \cite{Suy}, and \cite{Tsallis}.

\medskip

We consider the thermodynamic semiring as in Definition \ref{thermosemidef}
with the information measure $S$ given by the Tsallis entropy $S={\rm Ts}_\alpha$.

In this case the failure of the associativity condition for the semiring with the 
$\oplus_{S,T}$ operation is measured by comparing the expressions
\[  {\rm Ts}_{\alpha} (p_1) + (1 - p_1) {\rm Ts}_{\alpha} ( \frac{p_2}{1 -
p_1}) = \] 
\[ \frac{1}{\alpha - 1} (1 - p_1^{\alpha} - (1 - p_1)^{\alpha} +
\frac{p_2^{\alpha}}{(1 - p_1^{})^{\alpha - 1}} + \frac{(1 - p_1 -
p_2)^{\alpha}}{(1 - p_1)^{\alpha - 1}}) \]
and
\[ {\rm Ts}_{\alpha} (p_1 + p_2) + (p_1 + p_2) {\rm Ts}_{\alpha} (
\frac{p_1}{p_1 + p_2}) = \] 
\[ \frac{1}{\alpha - 1} (1 - (p_1 + p_2)^{\alpha} - (1 -
p_1 - p_2)^{\alpha} + \frac{p_1^{\alpha}}{(p_1 + p_2)^{\alpha - 1}} +
\frac{p_2^{\alpha}}{(p_1 + p_2)^{\alpha - 1}}) .\]

However, an interesting feature of the Tsallis entropy is that the
associativity of the thermodynamic semiring can be restored by
a deformation of the operation $\oplus_{S,T}$, depending on the
deformation parameter $\alpha$ which makes sense in the previous thermodynamic context, so that the Tsallis entropy
becomes the unique function that makes the resulting $\oplus_{S,T,\alpha}$ 
both commutative and associative.

\subsection{A Witt construction for Tsallis entropy}\label{WittTsallisSec}

We show here how to deform the thermodynamic semiring structure
in a one-parameter family $\oplus_{S,T, \alpha}$
for which $S={\rm Ts}_\alpha$ is the only entropy function that satisfies
the associativity constraint, along with the commutativity and unity axioms.

\smallskip

We consider here a slightly more general form of the Tsallis entropy,
as the non-associative information measure the
Tsallis entropy \cite{GeTsa}, defined by
\begin{equation}\label{phiTsallis}
{\rm Ts}_\alpha (p) = \frac{1}{\phi (\alpha)} (p^\alpha + (1 - p)^\alpha - 1), 
\end{equation}
where $\alpha \in \R$ is a parameter and $\phi$ is a continuous function
such that $\phi (\alpha) (1 - \alpha) > 0$, whenever $\alpha \neq 1$, with 
\[ \lim_{\alpha \rightarrow 1}
\phi (\alpha) = 0, \] 
and such that there exists $0 \leqslant a < 1 < b$ with the
property that $\phi$ is
differentiable on $(a, 1) \cup (1, b)$, and 
\[ \lim_{\alpha \rightarrow 1} \frac{d
\phi (\alpha)}{d \alpha} < 0. \] 

Note that this implies that the Tsallis entropy
reproduces the Shannon entropy in the $\alpha \rightarrow 1$ limit. 
A typical choice for the normalization is $\phi (\alpha) = 1 - \alpha$, which
reproduces the form \eqref{Tsallis}.

Here we work with the more general form \eqref{phiTsallis}, as we
will be able to ensure uniqueness only up to a general $\phi$ satisfying the above
requirements. 

We find that the Tsallis entropy fits nicely into the
context of Witt rings with the following two results.

\begin{thm}\label{uniqueTsallis}
  The Tsallis entropy in the form 
  \eqref{phiTsallis} is the unique entropy function that is commutative, has
  the identity property, and satisfies the $\alpha$-associativity condition
  \begin{equation}\label{qassoc}
   S (p_1) + (1 - p_1)^\alpha S ( \frac{p_2}{1 - p_1}) = S (p_1 + p_2) + (p_1 +
     p_2)^\alpha S ( \frac{p_1}{p_1 + p_2}) . 
  \end{equation}   
\end{thm}

\proof We assume a priori that $-S$ is concave and continuous. Therefore, $-S$
has a unique maximum, which is positive when $S$ is non-trivial, since $S (0)
= 0$. Moreover, $S$ is symmetric, so this maximum must occur at $p = 1 / 2$. 
$S$ also has the identity property and the $\alpha$-associativity, so by Suyari \cite{Suy}
and Furuichi \cite{Furu}, this implies $S = {\rm Ts}_\alpha$, for some $\phi (\alpha)$ 
satisfying the above properties. The converse follows from direct
application of the arguments given in \cite{Furu} and \cite{Suy}
and is easily verified. 
\endproof

\smallskip

The $\alpha$-associativity condition as one of the characterizing properties
for the Tsallis entropy was also discussed in \cite{CuTsa}.

\smallskip

We can interpret this $\alpha$-associativity as an associativity of an
$\alpha$-deformed Witt operation as follows. Fix some $\alpha$ and consider
\begin{equation}\label{xplusqy}
 x \oplus_{S,T,\alpha} y = \sum_{s \in I} e^{T S (s)} x^{s^\alpha} y^{(1 -
   s)^\alpha} . 
\end{equation}

We then have the following characterization of associativity.
   
\begin{thm}
  For $\alpha \neq 0$, the operation $\oplus_{S,T,\alpha}$ is associative 
  if and only if $S$ is $\alpha$-associative, as in \eqref{qassoc}.
\end{thm}

\proof We find that this operation is associative if and only if
\[ \sum_{s, r \in I} e^{T (S (s r) + (1 - s r)^\alpha S
   ( \frac{s (1 - r)}{1 - s r}))} x^{(sr)^\alpha}
   y^{(s (1 - r))^\alpha} z^{(1 - r)^\alpha} \]
\[ = \sum_{s, r \in I} e^{T (S (s) + s^\alpha S (r))}
   x^{(s r)^\alpha} y^{(s (1 - r))^\alpha} z^{(1 - r)^\alpha} . \]
   
We make the same subsitution as earlier, setting $p_1 = s r$, $p_2 = s
(1 - r)$, $p_3 = 1 - r$. Then the above condition becomes
\[ \sum_{p_1 + p_2 + p_3 = 1} e^{T (S (p_1) + (1 - p_1)^\alpha S ( \frac{p_2}{1 -
   p_1}))} x^{p_1^\alpha} y^{p_2^\alpha} z^{p_3^\alpha} \]
\[ = \sum_{p_1 + p_2 + p_3 = 1} e^{T (S (p_1 + p_2) + (p_1 + p_2)^\alpha S (
   \frac{p_1}{p_1 + p_2}))} x^{p_1^\alpha} y^{p_2^\alpha} z^{p_3^\alpha} . \]
   
When $\alpha \neq 0$, the map $a \mapsto a^\alpha$ is invertible and convex/concave, and the
above is a composition of this map with several Legendre transformations, so
we can invert this composition to obtain
\[ S (p_1) + (1 - p_1)^\alpha S ( \frac{p_2}{1 - p_1}) = S (p_1 + p_2) + (p_1 +
   p_2)^\alpha S ( \frac{p_1}{p_1 + p_2}), \]
which is exactly the $\alpha$-associativity condition.
\endproof

\medskip

It is worth pointing out at this point that in the above deformed Witt construction, we have replaced the energy functional
\[ U = \sum p_i E_i \]
with
\[ U_{\alpha} = \sum p_{i}^q E_i, \]
according to our interpretation in \ref{statmechSec}. In the setting of nonextensive statistical mechanics built upon the Tsallis entropy, this latter expression is exactly the energy functional used. Therefore, the deformed Witt addition is again naturally interpreted as a free energy, now in the more general $q$-deformed thermodynamics.

\section{The Kullback--Leibler divergence}\label{KLsec}

We now discuss another class of thermodynamic semirings in which
both the associativity and the commutativity properties fail, but in which we
can encode entropy functions varying over
some underlying space or manifold. In particular, we will connect the
thermodynamic semirings we consider in this section to the general
point of view of {\em information geometry}, as developed in \cite{AmaNa}, 
\cite{IkTaAma}.

\smallskip

The Kullback--Liebler divergence \cite{Kull}, \cite{KullLei} is a measure
of {\em relative entropy}, measured by the average logarithmic difference 
between two probability distributions $p$ and $q$. Since the averaging is done with
respect to one of the probability distributions, the KL divergence is not
a symmetric function of $p$ and $q$. 

\smallskip

More precisely, the KL divergence of two binary probability distributions
$p$ and $q$ is defined as
\begin{equation}\label{KLdiv}
{\rm KL} (p ; q) = p \log \frac{p}{q} 
+ (1 - p) \log \frac{1 - p}{1 - q}.
\end{equation}

The negative of the Kullback--Liebler divergence reduces to the Shannon entropy 
(up to a constant) in the case where $q$ is a uniform distribution. 
It is also called the {\em information gain}, in the sense that it measures 
the probability law $p$ relative to a given input or reference probability $q$.

\smallskip

We are especially interested here in considering the case where the
probability distribution $q$ depends on an underlying space of parameter,
continuously or smoothly. Mostly, we will be considering the following
two cases.

\begin{defn}\label{infomanifold}
A smooth univariate binary statistical $n$-manifold $\cQ$ is a set
of binary probability distributions $\cQ=(q(\eta))$ smoothly parametrized by 
$\eta\in \R^n$. 

A topological univariate binary statistical $n$-space $\cQ$ is a
set of binary probability distributions $\cQ=(q(\eta))$ continuously
parameterized by $\eta\in \Xi$, with $\Xi$ a compact Hausdorff
topological space.
\end{defn}

The first case leads to the setting of information geometry \cite{AmaNa}, 
\cite{IkTaAma}, while the second case is more suitable for treating 
multifractal systems \cite{BeSch}. 

\smallskip

We then consider thermodynamic semiring in the more general form
of Definition \ref{funthermosemi}. Let $\cX$ be either a compact subset
of $\R^n$ in the case of a smooth univariate binary statistical manifold
or a closed subset of a compact Hausdorff space $\Xi$ in the
topological case of Definition \ref{infomanifold}. We consider
the space of continuous functions $\cR=C(\cX,R)$, where the semiring
$K$ is either $\R^{\min,+}\cup \{ \infty \}$ or $\R^{\max,\ast}_{\geq 0}$, 
or in the smooth case we take $\cR=C^\infty(\cX,K)$. 

\smallskip

Give $q=q(\eta)$ in $\cQ$, we can endow the space $\cR$ of functions 
with a thermodynamic semiring structure as in Definition \ref{funthermosemi},
where the deformed addition operation is given by
\begin{equation}\label{plusKL}
x(\eta) \oplus_{{\rm KL}_\eta,\rho} y(\eta) = \sum_{p \in \Q \cap [0, 1]} \rho^{-{\rm KL} (p ;
q(\eta))} x(\eta)^p y(\eta)^{1 - p}, 
\end{equation}
where $\rho$ is the parameter of the deformation. Note we use the negative of the KL divergence because we are interested in it as a measure of relative entropy, rather than relative information, concepts often conceptually distinct but always related by a minus sign.

\smallskip

In the case when $q(\eta) \equiv 1 / 2$ is uniform for all $\eta$, 
we obtain back the original case with the Shannon entropy up to a shift factor
\[  x \oplus_{{\rm KL}_\eta,\rho} y |_{q(\eta)\equiv 1/2} 
= \max_p (- \rho (p \log \frac{p}{1 / 2} + (1 - p) \log \frac{1 -
p}{1 / 2}) + p x + (1 - p) y) \]
\begin{equation}\label{Shshift}
  = \max_p ( \rho {\rm Sh} (p) + p x + (1 - p) y) + \rho \log 2. 
\end{equation}  

\smallskip

We note that we can calculate this operation explicitly over $\R^{\min,+}\cup \{ \infty \}$ and $\R^{\max,\ast}_{\geq 0}$. We obtain the following result, by arguing as in Proposition \ref{propShsum}.

\begin{prop}\label{KLexpl}
We have the following expression over $\R^{\min,+}\cup \{ \infty \}$

\[x\oplus_{{\rm KL}} y = - T \log (e^{-\frac{x}{q T}}+e^{-\frac{y}{(1-q) T}})\]

and the following expression over $\R^{\max,\ast}_{\geq 0}$.

\[x\oplus_{{\rm KL}} y = ((\frac{x}{q})^{1/T}+(\frac{y}{1-q})^{1/T})^T\]
\end{prop}

The first observation then is that the additive structures \eqref{plusKL} are
in general not commutative. 

\begin{prop}\label{KLncomm}
The thermodynamic semiring structure 
\begin{equation}\label{KLsemi}
 x\oplus_{{\rm KL}} y = \sum_{p \in \Q \cap [0, 1]} \rho^{-{\rm KL} (p ;
q) } x^p y^{1 - p}
\end{equation}
is commutative if and only if $q=1/2$. The lack of commutativity
is measured by the transformation $q\mapsto 1-q$.
\end{prop}

\proof  This is immediate from the previous calculation, but we perform the proof over general $K$. We find
\[  {\rm KL} (1 - p ; q) = (1 - p) \log \frac{1 - p}{q} +
p \log \frac{p}{1 - q}. \]
This is related to ${\rm KL}(p;q)$ by the transformation
$q\mapsto 1-q$.
Thus, ${\rm KL} (p ; q) = {\rm KL} (1 - p ; q)$ when $\log \frac{1 - q}{q} =
0$, that is, when $q = 1 / 2$. This is exactly when
the Shannon entropy case is reproduced, so the only case when the
addition \eqref{KLsemi} based on the Kullback--Liebler divergence
is commutative is when it agrees with the Shannon entropy up to a shift factor.
\endproof

For the associativity condition we find the following result.

\begin{prop}\label{KLassoc}
The lack of associativity of the thermodynamic semiring \eqref{KLsemi}
is measured by the transformation $(p_1,p_2,p_3;q)\mapsto (p_3,p_2,p_1;1-q)$.
\end{prop}

\proof Again we proceed over general $K$. We have
\[ {\rm KL} (p_1 ; q) + (1 - p_1) {\rm KL} ( \frac{p_2}{1 - p_1} ; q) \]
\[  = p_1 \log \frac{p_1}{q} + (1 - p_1) \log \frac{1 - p_1}{1 - q} \]
\[ + p_2 \log \frac{p_2}{(1-p_1) q} 
+ (1-p_1-p_2)  \log \frac{1 - p_1 - p_2}{(1
- p_1) (1 - q)} \]
\[ =  p_1 \log \frac{p_1}{q} + (1 - p_1) \log \frac{1 - p_1}{1 - q}  \]
\[ + p_3 \log \frac{p_3}{1-q} + p_2 \log \frac{p_2}{q} - (1-p_1) \log (1-p_1), \]
while
\[ {\rm KL} (p_1 + p_2 ; q) + (p_1 + p_2) {\rm KL} ( \frac{p_1}{p_1 + p_2} ;
q) \]
\[ = (p_1 + p_2) \log \frac{p_{_1} + p_2}{q} + (1 - p_1 - p_2) \log \frac{1 - p_1 -
p_2}{1 - q} \] \[+ (p_1 + p_2) \frac{p_1}{p_1 + p_2} \log \frac{p_1}{(p_1 + p_2) q}
+ (p_1 + p_2) \frac{p_2}{p_1 + p_2} \log \frac{p_2}{(p_1 + p_2) (1 - q)} \]
\[ = (p_1 + p_2) \log \frac{p_1 + p_2}{q} + (1 - p_1 - p_2) \log \frac{1 - p_1 -
p_2}{1 - q} \] \[+ p_1 \log \frac{p_1}{q} + p_2 \log \frac{p_2}{1 - q} - (p_1 + p_2)
\log (p_1 + p_2) \]
\[ = p_3 \log\frac{p_3}{1-q} + (1-p_3) \log\frac{1-p_3}{q}  \]
\[  + p_1 \log \frac{p_1}{q} + p_2 \log \frac{p_2}{1-q} - (1-p_3) \log (1-p_3). \]
These are related by the 
transformation $(p_1,p_2,p_3;q)\mapsto (p_3,p_2,p_1;1-q)$. 
\endproof

Notice that, because of the presence of the shift in \eqref{Shshift} with
respect to the Shannon entropy, in the case $q=1/2$ we find
\[ {\rm KL} (p_1 ; \frac{1}{2}) + (1 - p_1) {\rm KL} ( \frac{p_2}{1 - p_1} ; \frac{1}{2}) = \]
\[ p_1 \log p_1 + p_2 \log p_2 + p_3 \log p_3 + \log 2 + (1-p_1) \log 2 \]
while
\[  {\rm KL} (p_1 + p_2 ; \frac{1}{2}) + (p_1 + p_2) {\rm KL} ( \frac{p_1}{p_1 + p_2} ;
\frac{1}{2})  =\]
\[ p_1 \log p_1 + p_2 \log p_2 + p_3 \log p_3 + \log 2 + (1-p_3) \log 2 . \]

Thus, associativity is not automatically obtained in the uniform distribution case, but instead we have associativity up to a shift.

\smallskip

By Proposition \ref{KLassoc} we see that, 
in the case of a thermodynamic semiring $\cR=C(\cX,K)$
or $\cR=C^\infty(\cX,K)$, for a topological or smooth univariate binary
statistical space, if one can find an involution $\alpha: \cX \to \cX$ of the
parameter space such that $q(\alpha(\eta))=1-q(\eta)$, then one can
consider the transformation $x(\eta) \mapsto x(\alpha(\eta))$ and one
finds 
\[ x(\eta)\oplus_{{\rm KL}_{q(\eta)}} y(\eta)=y(\alpha(\eta))\oplus_{{\rm KL}_{q(\alpha(\eta))}} x(\alpha(
\eta)) . \]
Moreover, the morphism 
$$ A: (x(\eta),y(\eta),z(\eta))\mapsto (z(\alpha(\eta)),y(\alpha(\eta)),x(\alpha(\eta))) $$
measures the lack of associativity, by making the diagram commute,
$$\xymatrix{ \cR\otimes \cR \otimes \cR    \ar@{->}[rr]^{A} \ar@{->}[d]_{\oplus_{{\rm KL}}\otimes 1} & &
\cR\otimes \cR \otimes \cR \ar@{->}[d]^{1\otimes \oplus_{{\rm KL}}}  \\ 
\cR\otimes \cR \ar@{->}[r]^{\oplus_{{\rm KL}}} & \cR & \cR\otimes \cR \ar@{->}[l]_{\oplus_{{\rm KL}}} } $$

\smallskip

\subsection{Applications to multifractal systems}

Consider the case of a Cantor set $\cX$ identified, through
its symbolic dynamics interpretation, as the one sided full
shift space $\Sigma_2^+$ on the alphabet $\{ 0,1 \}$, see \S 1.3
of \cite{PeCli}.

\smallskip

For $\eta \in \cX$, let $a_n(\eta)$ denote the number of $1$'s that
appear in the first $n$ digits $\eta_1, \ldots, \eta_n$ of $\eta$. We set
\begin{equation}\label{qxdigits}
q(\eta) = \lim_{n\to \infty} \frac{a_n(\eta)}{n},
\end{equation}
if this limit exists.  We denote by $\cY \subset \cX$ the set of
points for which the limit \eqref{qxdigits} exists.

\smallskip

The limit \eqref{qxdigits} determines several important dynamical
properties related to the fractal geometry of $\cX$. For example, 
suppose that $\cX$ is a uniform Cantor set obtained from a contraction 
map $f$ with contraction ratio $\lambda$, endowed with 
a Bernoulli measure $\mu_p$ for a given $0<p<1$, 
defined by assigning measure 
$$ \mu_p(\cX(w_1,\ldots,w_n))=p^{a_n(w)} (1-p)^{n- a_n(w)} $$
to the cylinder sets 
$$ \cX(w_1,\ldots,w_n) =\{ \eta\in \cX \,|\, \eta_i =w_i, \, i=1,\ldots,n \}. $$
Then, the local dimension of $\cX$ at a point $\eta\in \cY$ is given by (\S 4.17 of \cite{PeCli})
$$ d_{\mu_p}(\eta) = \frac{q(\eta) \log p + (1-q(\eta)) \log (1-p)}{\log \lambda} $$
while the local entropy of the map $f$ is given by  (\S 4.18 of \cite{PeCli})
$$ h_{\mu_p,f}(\eta) = q(\eta) \log p + (1-q(\eta)) \log (1-p). $$
For a non-uniform Cantor set $\cX$ with two contraction ratios $\lambda_1$
and $\lambda_2$ on the two intervals, the Lyapunov exponent of $f$ is
given by (\S 4.20 of \cite{PeCli})
$$ \lambda_f(\eta) = q(\eta) \log \lambda_1 + (1-q(\eta)) \log \lambda_2 . $$

\smallskip

One knows that, given a Bernoulli measure $\mu_p$ on the Cantor set $\cX$,
there is a set $\cZ\subset \cX$ of full measure $\mu_p(\cZ)=1$, for which
$q(\eta)=p$ (Proposition 4.5 of \cite{PeCli}). The choice of the uniform measure
$\mu_{1/2}$ yields a full measure subset $\cZ_{1/2}$ on which the limit $q(\eta)=1/2$
is the uniform distribution (the fair coin case). In general one can
stratify the set $\cY\subset \cX$ into level sets of $q(\eta)$. This provides a
decomposition of the Cantor set as a multifractal.  

\smallskip

Looking at this setting from the point of view of thermodynamic semirings
suggests considering the set of functions $C(\cY,K)$ endowed with the
pointwise operation $\oplus_{{\rm KL}_{q(\eta)},T}$, with the Kullback--Leibler 
divergence ${\rm KL}(p;q(\eta))$, for $q(\eta)$ defined as in \eqref{qxdigits}.
Then we see that, without the need to choose a measure on $\cX$, the
algebraic properties of the thermodynamic semiring automatically select
the ``fair coin subfractal" $\cZ_{1/2}$.

\begin{prop}\label{multifrprop}
For $\cZ \subset \cY$, the semiring $C(\cZ,K)$, with the operation 
$\oplus_{{\rm KL}_{q(\eta)},T}$, for $q(\eta)$ as in \eqref{qxdigits}, is 
commutative if and only if $\cZ\subset \cZ_{1/2}$ is a ``fair coin" subset.
\end{prop}

\proof It follows immediately from Proposition \ref{KLncomm}.
\endproof

Moreover, we can see geometrically the involution that
measures the lack of commutativity as in Proposition \ref{KLncomm}
and the lack of associativity as in Proposition  \ref{KLassoc}.

\begin{prop}\label{assocCantor}
The homeomorphism $\gamma: \cX \to \cX$ given by the involution 
that exchanges $0\leftrightarrow 1$ in the digits of $\eta$ in the
shift space $\Sigma_2^+$ implements the involution $q(\eta)\mapsto
1-q(\eta)$ that measures the lack of commutativity and that, together
with the involution $(p_1,p_2,p_3)\mapsto (p_3,p_2,p_1)$ also
measures the lack of associativity. Thus, the morphism $x(\eta)\mapsto
x(\gamma(\eta))$ restores commutativity, in the sense that 
$$ x(\eta)\oplus_{{\rm KL}_{q(\eta)}} y(\eta)=
y(\gamma(\eta))\oplus_{{\rm KL}_{q(\gamma(\eta)}} x(\gamma(\eta)),$$ 
while
$A: \cR\otimes \cR \otimes \cR\to \cR\otimes \cR \otimes \cR$
given by $$A(x(\eta),y(\eta),z(\eta))=(z(\gamma(\eta)),y(\gamma(\eta)),x(\gamma(\eta)))$$
restores associativity, making the diagram commute
$$\xymatrix{ \cR\otimes \cR \otimes \cR    \ar@{->}[rr]^{A} \ar@{->}[d]_{\oplus_{{\rm KL}}\otimes 1} & &
\cR\otimes \cR \otimes \cR \ar@{->}[d]^{1\otimes \oplus_{{\rm KL}}}  \\ 
\cR\otimes \cR \ar@{->}[r]^{\oplus_{{\rm KL}}} & \cR & \cR\otimes \cR \ar@{->}[l]_{\oplus_{{\rm KL}}} } $$
\end{prop}

\proof This follows immediately from Proposition \ref{KLncomm} and
Proposition  \ref{KLassoc}, by observing that the $q(\eta)$ defined as in
 \eqref{qxdigits} satisfies $q(\gamma(\eta))=1-q(\eta)$, since $a_n(\gamma(\eta))=n -a_n(\eta)$
 for all $\eta \in \cX$.
\endproof 

\medskip
\subsection{Multivariate binary statistical manifolds}

We see that in the univariate case, the extremal $p$ value is the unique probability distribution minimizing the KL-divergence to $q$ subject to the soft constraint coming from the energy functional $p x + (1-p) y$, see \S \ref{statmechSec}. This is important because minimizing the KL
divergence is maximizing likelihood, and this plays an important role
in marginal estimation, belief propagation, mutual information calculation, 
see \cite{IkTaAma} and \cite{AmaNa}.
\smallskip

A more interesting case is that of multivariate statistical manifolds.
To maintain the same features as in the univariate case, we will find that a hyperring structure is most natural. See \cite{Vi}, \cite{CC} for an introduction and relevent facts of hyperstructures. We first note the following fact.

\begin{prop}\label{KLmulti}
  If $p$ and $q$ are two distributions, we denote by $p_i$ and $q_i$ their
  $i$-th marginal distribution. Then ${\rm KL} (p ; q) = \sum_i {\rm KL}
  (p_i ; q_i)$.
\end{prop}

\proof We have
\begin{align*} 
{\rm KL} (p ; q) = & p_1 \cdots p_n \log \frac{p_1 \cdots p_n}{q_1 \cdots q_n} \\
+&  (1 - p_1) p_2 \cdots p_n \log \frac{(1 - p_1) p_2 \cdots p_n}{(1 - q_1) q_2
\cdots q_n} \\ 
+ \cdots + & (1 - p_1) \cdots (1 - p_n) \log \frac{(1 - p_1) \cdots
(1 - p_n)}{(1 - q_1) \cdots (1 - q_n)}  
\end{align*}
\[
=  p_1 \cdots p_n (\log \frac{p_1}{q_1} + \cdots + \log \frac{p_n}{q_n}) + \cdots
+ (1 - p_1) \cdots (1 - p_n) (\log \frac{1 - p_1}{1 - q_1} + \cdots + \log
\frac{1 - p_n}{1 - q_n}) 
\]
\begin{align*}
 = & p_1 \log \frac{p_1}{q_1} (p_2 \cdots p_n + (1 - p_2) \cdots p_n + \cdots) \\
 + & (1 - p_1) \log \frac{1 - p_1}{1 - q_1} (p_2 \cdots p_n + \cdots) + \cdots 
 \end{align*}
\[ = p_1 \log \frac{p_1}{q_1} ((1 + p_2 - p_2) (p_3 \cdots p_n + \cdots)) + \cdots \]
\[ = p_1 \log \frac{p_1}{q_1} + (1 - p_1) \log \frac{1 - p_1}{1 - q_1} + \cdots +
(1 - p_n) \log \frac{1 - p_n}{1 - q_n} = \sum_i {\rm KL} (p_i ; q_i). \]
\endproof

Thus, if we can ensure that the sum of the KL divergences of the marginal
distributions is minimized, then the total KL divergence will be minimized.

\subsection{Product of semirings and hyperfield structure}

We proceed by taking the semiring
$$ \cR=C(\{ 1,\ldots, n\}, K)=K^{\otimes n}. $$ It is tempting to define the operations on $\cR$ coordinate-wise, however, since we want to consider an $n$-ary probability distribution and not $n$ binary probability distributions, there should be some dependence between coordinates that takes advantage of the previous proposition. In short, we would like to put an ordering on $\cR$ that ensures the trace
\[(x_1,...,x_n) \rightarrow x_1+...+x_n \in K \]
is maximized. This ordering does not uniquely determine a maximum between two tuples. We thus forsake well-definedness of the addition on $K$ and define $(x_1,...,x_n)+(y_1,...,y_n)$ to be the set of tuples $(z_1,...,z_n)$ with $z_i = x_i$ or $y_i$ that maximize $z_1+...+z_n$ in the ordering on $K$. This, together with coordinate-wise multiplication defines a characteristic one hyperfield structure on $\cR$. We then define the Witt operation for some information measures $S_1, ..., S_n$ over $K=\R^{\min,+}\cup \{ \infty \}$ by
\[x\oplus_{S_1,...,S_n}y = \min_{p_1,...,p_n} (p_1 x_1+(1-p_1) y_1 - T S_1 (p_1),...,p_n x_n+(1-p_n) y_n - T S_n (p_n)),\]
where $x=(x_1,...,x_n), y=(y_1,...,y_n)$, now we consider the $p_i$ as marginal probabilities, and the $\min$ operation is the multivalued hyperring addition. When each $S_i$ is the KL-divergence from some $q_i$, by the previous proposition, the results of this operation are exactly the distributions with marginal probabilities $(p_1,...,p_n)$ minimizing the KL-divergence to the marginal probabilities $(q_1,...,q_n)$ subject to the soft constraint coming from the energy functional
\[U = \sum p_i x_i + (1-p_i) y_i.\]
The lack of well-definedness of this addition can be interpreted in the thermodynamic context as the non-uniqueness of equilibria, via the existence of meta-equilibrium states. Indeed, when the $q_i$ describe a uniform distribution, we find that this addition is in fact well-defined.
\smallskip

Note that these hyperfields are slightly different from those considered in \cite{Vi}. However, just as taking $T \leftarrow 0$ for the Shannon entropy semiring reproduces the ``dequantized" tropical semiring, we can take $T \leftarrow 0$ for the KL divergence semiring to get a ``dequantized" tropical hyperfield. This reproduces the underformed addition defined on $K$ above. Note that this is not the same as Oleg Viro's tropical hyperfield discussed in \cite{Vi}.

\smallskip

We can encode more information about a space in
the ring deformation by restricting the marginal probabilities we
sum over, in particular we can restrict the minimizing process to certain submanifolds of our probability manifold
such as the e-flat or m-flat manifolds typically considered in \cite{AmaNa}, since the KL-divergence is related to the Fisher information matrix defining the Riemannian structure.
See also the comments in \S \ref{infogeomsec} below.

\section{The successor in thermodynamic semirings}\label{assocSec}

\smallskip

Given a thermodynamic semiring in the sense of Definition \ref{thermosemidef},
we let
\begin{equation}\label{lambdafirst}
\lambda (x,T) = x \oplus_{S} 0 \equiv \min_p (p x - T S (p)). 
\end{equation}
Then $\lambda : K \times \R \rightarrow K$ is the Legendre transform of $T S : [0, 1]
\rightarrow \R$. If we assume that $S$ has a unique maximum, then we
can invert the Legendre transform, so that
\[  T S (p) = \min_x (p x - \lambda (x,T)). \]
Therefore, when $S$ is concave/convex, we can recover it from the
semiring. We call $\lambda$ the successor function, since $0$ is the multiplicative identity, and over general $K$ we can write $\lambda (x,T) = x \oplus_S 1$. When multiplication distributes over addition, we can write
\[x\oplus_S y = \lambda (x-y,T)+y.\]
We will tend to suppress the $T$ dependence of $\lambda$. Each of the algebraic properties of $S$ and $K$ translate into the language of $\lambda$.

\smallskip

\begin{prop}\label{lambdaSaxioms}
The entropy function $S$ has the following properties.
\begin{enumerate}
\item It satisfies the commutativity axiom  
$S (p) = S (1 - p)$ (hence $\oplus_{S,T}$ is commutative) if and only if 
 \begin{equation}\label{commutlambda} 
   \lambda (x) - \lambda (-x) = x.
 \end{equation}  
\item It satisfies the left identity axiom $S (0) = 0$ (hence $\oplus_{S}$ has left identity $\infty$) 
if and only if $\lambda (x) \leqslant 0$
and $\lim_{x \rightarrow \infty}  \lambda (x) = 0$.
\item It satisfies the right identity axiom $S (1) = 0$ (hence $\oplus_{S}$ has left identity $\infty$) 
if and only if $\lambda (x) \leqslant x$ 
and $ \lambda (x) \sim x$, as $x \rightarrow - \infty$.
\item It satisfies the associativity constraint making  $\oplus_{S}$ associative 
iff
$$ \lambda (x - \lambda (y)) + \lambda (y) =
  \lambda (\lambda (x - y) + y). $$
\end{enumerate}
\end{prop}

\proof  Facts (1) and (4) are immediate from the definition. The properties (2) 
and (3) arise from the fact that $\lambda$ should be continuous at 
$\infty$ and $- \infty$. We then read $\infty \oplus_{S} x$ and $x \oplus_{S} \infty$ as 
$\lim_{y \rightarrow \infty} y \oplus_{S} x$ and $\lim_{y \rightarrow \infty} x \oplus_{S} y$, 
respectively. Each of these should equal $x$, and in terms of $\lambda$ 
we see that $\lim_{y \rightarrow \infty} y \oplus_{S} x = \lim_{y \rightarrow \infty} \lambda (y - x) +
x$ and $\lim_{y \rightarrow \infty} x \oplus_{S,T} y = \lim_{y \rightarrow \infty}
\lambda (x - y) + y$, thus proving (2) and (3).
\endproof

\smallskip

In the case of the Shannon entropy $S={\rm Sh}$ and KL-divergence $S=-{\rm KL}(p;q)$, we have the following
forms for the successor function.

\begin{prop}\label{Shannonlambda}
For Shannon entropy,
\begin{equation}\label{lambdaSh}  
 \lambda^{{\rm Sh}} (x,T) = - T \log (1 + e^{- x / T}),
\end{equation} 
over $\R^{\min, +}  \cup \{\infty\}$, and
\begin{equation}\label{lambdaSh2}  
\lambda^{{\rm Sh}} (x,T) = (1 + x^{1 / T})^T 
\end{equation}
over $\R^{\max, +}_{\geqslant 0}$. For the KL-divergence,
\begin{equation}\label{lambdaKL}
 \lambda^{{\rm KL}} (x,T) = - T \log (1 + e^{- x/q T}),
\end{equation}
over $\R^{\min, +}  \cup \{\infty\}$, and
\begin{equation}\label{lambdaKL}
 \lambda^{{\rm KL}} (x,T) = (1/(1-q)^{1/T}+(x/q)^{1/T})^T,
\end{equation}
over $\R^{\max, +}_{\geqslant 0}$.
\end{prop}

\proof  This follows directly from the
definition of $\lambda (x,T) = x \oplus_{S} 0$, and the isomorphism
$-\log$ relating the semirings 
$\R^{\max,*}_{\geq 0}$ and $\R^{\min,+}\cup \{ \infty \}$.
\endproof

\begin{figure}[h]
 \includegraphics{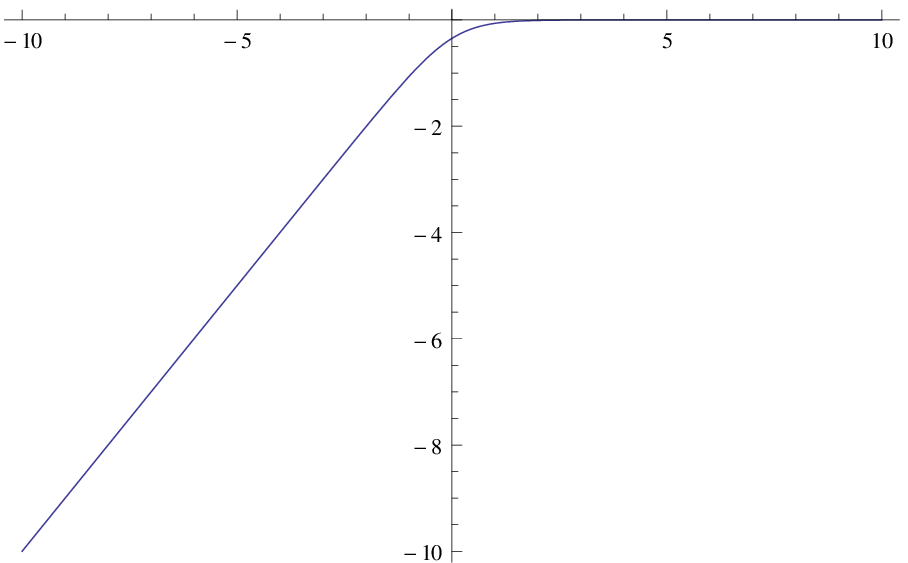}
 \caption{The successor function $\lambda^{{\rm Sh}}$ for $T=0.5$ 
 \label{figSh1}}
\end{figure}

\begin{figure}[h]
 \includegraphics{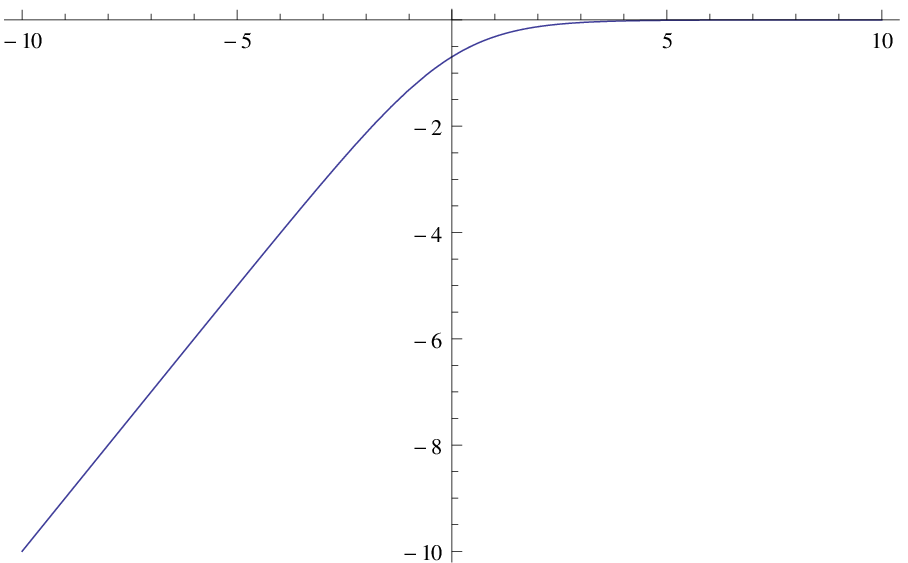}
 \caption{The successor function $\lambda^{{\rm Sh}}$ for $T=1$ 
 \label{figSh2}}
\end{figure}

\begin{figure}[h]
 \includegraphics{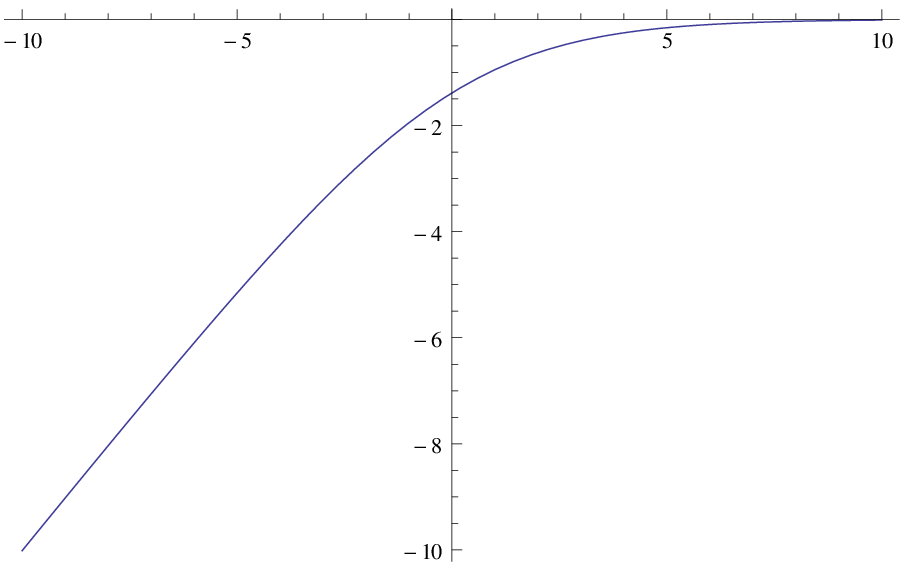}
 \caption{The successor function $\lambda^{{\rm Sh}}$ for $T=2$ 
 \label{figSh3}}
\end{figure}

Figures \ref{figSh1}, \ref{figSh2} and \ref{figSh3}
show examples of a plot of $\lambda^{{\rm Sh}}$ plotted 
versus $x$, for different values of $T$.

\subsection{Successor function for Tsallis entropy}

Consider now the case of the Tsallis entropy
\[ {\rm Ts}_{\alpha} (p) = \frac{1}{1 - \alpha} (p^{\alpha} + (1 - p)^{\alpha}
- 1). \]

\begin{prop}\label{TsallisAssoc}
The successor function $\lambda^{{\rm Ts}_\alpha} (x,T)$ for the Tsallis entropy
is given by
\begin{equation}\label{eqTsAssoc} 
\lambda^{{\rm Ts}_\alpha} (x, T) = \left\{ \begin{array}{lcl}
  0   &  &   | \frac{\alpha}{1 - \alpha} | < x / T  \\[3mm]  
  g (x)  &  &  - | \frac{\alpha T}{1 - \alpha} | < x / T < | \frac{\alpha}{1 - \alpha} |  \\[3mm]  
  x   &  &  x / T < - | \frac{\alpha T}{1 - \alpha} |
\end{array} \right.  
\end{equation}
where $g (x)$ is given by applying ${\rm Ts}$ to the inverse of its
derivative. 
\end{prop}

\proof We have
\[ \frac{\partial {\rm Ts}}{\partial p} = \frac{\alpha}{1 - \alpha} (p^{\alpha
- 1} - (1 - p)^{\alpha - 1}) . \]
We see the derivative of ${\rm Ts}_{\alpha}$ has range $[- | \frac{\alpha}{1
- \alpha} |, | \frac{\alpha}{1 - \alpha} |]$, so that we obtain \eqref{eqTsAssoc}.
\endproof

\begin{figure}[h]
 \includegraphics{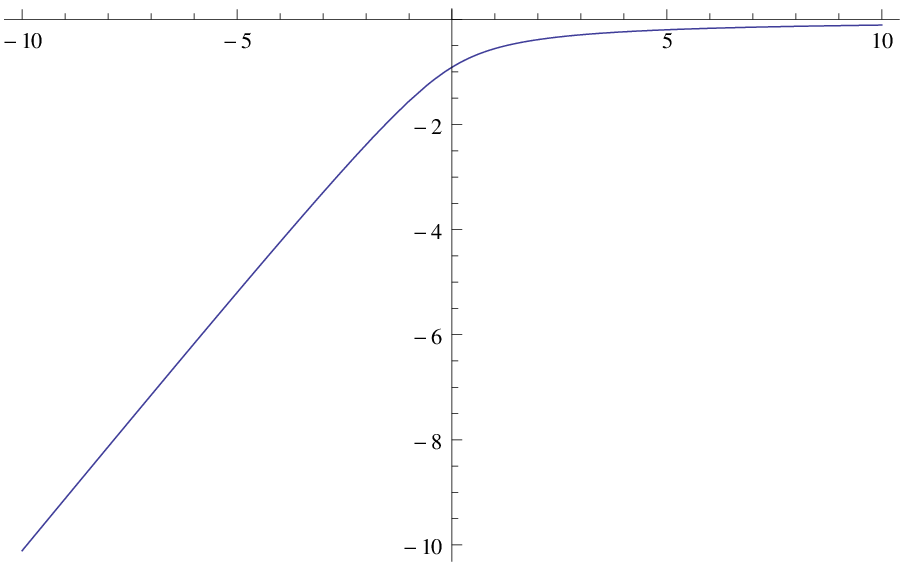}
 \caption{The successor function $\lambda^{{\rm Ts}}_{\alpha}$ for $\alpha=0.5$ and $T=1$ 
 \label{figTs1}}
\end{figure}

Figure \ref{figTs1} shows an example of a plot of $\lambda^{{\rm Ts}_\alpha}$ plotted 
versus $x$. In the limit $\alpha \rightarrow \infty$,  one has 
${\rm Ts}_{\alpha} =\chi_{(0, 1]}$, so indeed $\lambda^{{\rm Ts}_{\infty}} (x) = x \chi_{[-
\infty, 0)} (x)$ for finite temperature. When $\alpha < 0$,
${\rm Ts}_{\alpha}$ is convex, so $\lambda$ becomes concave in this region,
as expected.

\medskip
\subsection{Successor function for R\'enyi entropy}

We now consider again the R\'eyni entropy given by
\[ {\rm Ry}_{\alpha} (p) = \frac{1}{1 - \alpha} \log (p^{\alpha} + (1 -
p)^{\alpha}). \]
We have
\[ \frac{\partial {\rm Ry}}{\partial p} = \frac{\alpha}{1 - \alpha} (p^{\alpha
- 1} + (1 - p)^{\alpha - 1}) / (p^{\alpha} + (1 - p)^{\alpha}).  \]

This time, however, the derivative has range $\R$, so that we have both
$\lambda^{{\rm Ry}_{\alpha}} (x) < x$ and $\lambda^{{\rm Ry}_\alpha}
(x) < 0$.

\smallskip

\begin{figure}[h]
 \includegraphics{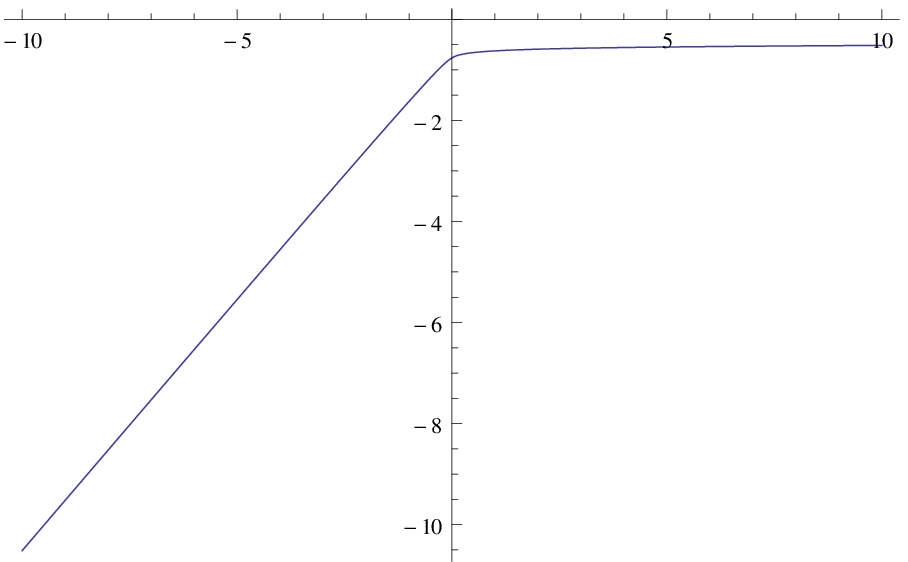}
 \caption{The successor function $\lambda^{{\rm Ry}_\alpha}$ for $\alpha=0.1$ and $T=1$ 
 \label{figRy1}}
\end{figure}

\begin{figure}[h]
 \includegraphics{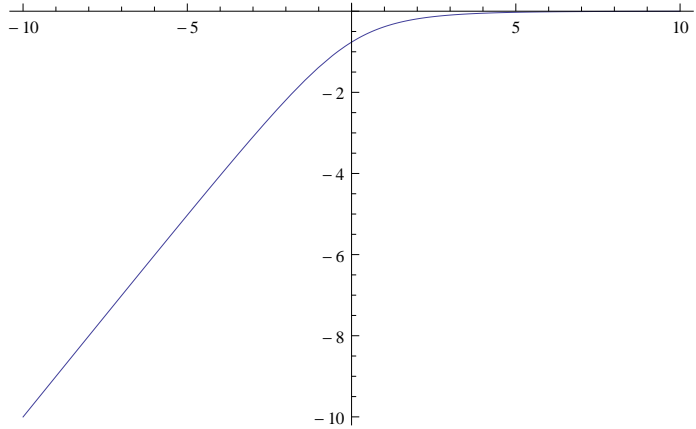}
 \caption{The successor function $\lambda^{{\rm Ry}_\alpha}$ for $\alpha=0.9$ and $T=1$ 
 \label{figRy2}}
\end{figure}

Figures \ref{figRy1} and \ref{figRy2}
show examples of a plot of $\lambda^{{\rm Ry}_\alpha}$ plotted 
versus $x$, for different values of $\alpha$.

One can see, by comparing these various graphs for the different entropy functions 
that increasing $T$ has the effect of smoothing the transition, 
while increasing $\alpha$ sharpens it. 

\smallskip
\subsection{Cumulants generating function}\label{cumulsec}

In this section we give a thermodynamical interpretation of the
successor function.

Recall that, for a random variable $X$, if $M_X(t)$ denotes the
generating function for the momenta of $X$,
$$ M_X(t)= \langle \exp(tX) \rangle = \sum_{m=0}^\infty \mu_m \frac{t^m}{m!}, $$
then the {\em cumulants} $\{ \kappa_n \}$ of $X$ are defined 
as the coefficients of the power series expansion of the function $\log M_X(t)$,
$$ \log M_X(t) = \sum_{n=0}^\infty \kappa_n \frac{t^n}{n!}. $$
The information contained in cumulants or momenta is equivalent, 
though cumulants have the advantage that they behave additively over
independent variables.

\smallskip

We then have the following result. We consider the case of an analytic 
$\lambda$, which is reasonable when attempting to gain a thermodynamic 
understanding, as the microscopic dynamics are usually assumed to be analytic.

\begin{prop}\label{thmcumulant}
Let $\lambda (x,T)$ be the successor function of a thermodynamic semiring $K$.
Assume that $\lambda (x,T)$ is analytic.
Then the function $-\lambda (x,T)/T$ is the cumulant generating function of the
probability distribution for the energy $E$, in the variable $- 1 / T = -
\beta$. Namely, if we write the $n$th cumulant as $\kappa_n=\langle E^n \rangle_c$, 
we have 
\begin{equation}\label{kappanE} 
(- 1)^{n + 1} \frac{\partial^n}{\partial \beta^n} (\beta \lambda (x,T)) =
   \langle E^n \rangle_c . 
\end{equation}   
\end{prop}

\proof
In thermodynamics, $Z(\beta)=\langle \exp(-\beta E) \rangle$ is the partition
function, where $\beta=1/T$ is the inverse temperature and $E$ is the energy.
The Helmholtz free energy is then given by
\begin{equation}\label{freeenergy}
  F = - T \log \langle \exp (- E / T) \rangle .
\end{equation}
Up to a factor of $-1/T$, the 
Helmholtz free energy is in fact the cumulant generating function for 
the random variable given by the energy $E$. As observed already
in \S \ref{statmechSec}, the Helmholtz free energy is the Legendre
transform of the entropy, and can therefore we identified, again
up to a factor of $-1/T$, with the function $\lambda (x,T)$.
\endproof

\smallskip
      
We can of course perform this proof without reference to the thermodynamics. That is to say: the Legendre transform structure of the whole ordeal is independent of the information measure we select as long as we select one which is concave and analytic.
\smallskip
In particular, from \eqref{kappanE} we have
\[ \lambda (x,T) - T \frac{\partial}{\partial T} \lambda (T,x) = \langle E
   \rangle = p_{{\rm eq}} x , \]
where $p_{{\rm eq}} = p_T (x)$ is equilibrium value of the mole fraction. We know that $\lambda (x,T) = 
\min_p (p x - T S (p)) = p_T (x) - T S (p_T (x))$. We see that $p_T (x)$ satisfies
\[ x / T = \frac{d}{d p} S (p_T (x)) , \]
so we can write $p_T (x) = p (x / T)$ and $\lambda (T,x) = \lambda (x / T)$.
Notice that this explains the effect that changing the temperature has on $\oplus_{S,T}$. 

\smallskip

{}From the definition, we calculate
\[ \frac{\partial}{\partial T} \lambda (x / T) = x \frac{\partial}{\partial T}
   p (x / T) - S (p (x / T)) - T \frac{\partial}{\partial T} p (x / T)
   \frac{d}{d p} S (p (x / T)) , \]
which, by the above property, is just $- S (p (x / T))$, proving the above
relation. Note this holds for arbitrary smooth, concave entropy functions.
Similarly, we calculate 
\[ \frac{\partial}{\partial x} \lambda (x / T) = x p (x / T), \] 
so that
\[ \lambda (x / T) = x \frac{\partial}{\partial x} \lambda (x / T) + T
   \frac{\partial}{\partial T} \lambda (x / T) . \]
This is a well-known property of the Legendre transform of smooth functions.

\section{Entropy Operad}\label{operadsec}

A categorical and operadic point of view on convex spaces and entropy functions
was recently proposed in \cite{BaFriLei}, \cite{BaFriLei2}, \cite{Fri}, \cite{Fri2}. Here
we will use a similar viewpoint to describe generalized associativity conditions on
thermodynamic semirings.

\smallskip

More precisely, we consider here the more general question of how binary 
(or more complicated) information measures can be built up to ones 
for $n$-ary random variables for any $n \geqslant 2$. 
This will give us some interesting correspondences
between the combinatorics of such ``guessing games'' and generalized
associativity conditions in an operad with $n$-ary operations defined over $K$
like $x_1 \oplus_S \cdots \oplus_S x_n$ with some choice of parenthesizing. In
this section, we will assume for simplicity that $K$ is $\R^{\min, +} \cup \{\infty\}$,
although, once again, this is only a notational convention chosen to elucidate
certain expressions. All the statements made here could be translated into the
greater generality for real characteristic one semifields.

\smallskip

Operads were first introduced in \cite{May2} in the theory of iterated loop spaces and
have since seen a broad range of applications in algebra, topology, and geometry.
We recall briefly some basic facts about operads that we will need later, see
\cite{May}.

\smallskip

An operad is a collection of objects from a
symmetric monoidal category $\mathcal{S}$ with product $\otimes$ and unit
object $\kappa$. In particular, for each $j \in \N$, we have an
object $\mathcal{C}(j)$, thought of as parameter objects for $j$-ary
operations, with actions by the symmetric group ${\rm Sym}_j$, thought of as
permuting inputs. We also have a unit map $\eta : \kappa \rightarrow
\mathcal{C}(1)$ and composition maps
\[ \gamma : \mathcal{C}(k) \otimes \mathcal{C}(j_1) \otimes \cdots \otimes
   \mathcal{C}(j_k) \rightarrow \mathcal{C}(j_1 + \cdots + j_k) \]
which are suitably associative, unital, and equivariant under the action of
${\rm Sym}_k$ such that if $\sigma \in {\rm Sym}_k$,
\[ \gamma (c_k \otimes c_{\sigma (j_1)} \otimes \cdots \otimes c_{\sigma
   (j_k)}) = \gamma (\sigma (c_k) \otimes c_{j_1} \otimes \cdots \otimes
   c_{j_k}) . \]
A $\mathcal{C}$-algebra $A$ is an object together with
${\rm Sym}_j$-equivariant maps
\[ \mathcal{C}(j) \otimes A^j \rightarrow A, \]
thought of as actions, which are suitably associative and unital. Here $A^j$
represents $A^{\otimes j}$ and $A^0 = \kappa$. An $A$-module $M$ is an object
together with ${\rm Sym}_{j - 1}$-equivariant maps
\[ \mathcal{C}(j) \otimes A^{j - 1} \otimes M \rightarrow M \]
which are also suitably associative and unital. Note that we are taking our objects
all from symmetric monoidal categories, so we do not need to distinguish where the
operad lives from where the algebras live, but we have not eliminated the
possibility of doing so. When we consider the entropy
operad, the $n$-ary operations of the operad will be parametrized by rooted
trees, while we will take algebras from the category of topological categories.

\smallskip

\subsection{Operads and entropy}

We first recall the recent construction of J. Baez, T. Fritz,
and T. Leinster, \cite{BaFriLei}, \cite{BaFriLei2} of an 
operadic formalism for entropy, which is especially relevant 
to our setting and nicely displays the basic machinery.

\smallskip

Using the set theorists' convention, we define natural numbers as $n = \{0, \ldots,
n - 1\}$. An ordered $n$-tuple will then be denoted as $(a_i)_{i \in n} = (a_0,
\ldots, a_{n - 1})$. Consider as our symmetric monoidal category the category
of topological categories, denoted ${\rm Cat} ({\rm Top})$, with $\otimes$
as the Cartesian product, and $\kappa$ as the one-point space. One can
construct an operad, $\mathcal{P}$, out of probability distributions on finite sets. 
For each $j$, we define $\mathcal{P}(j)$ as the set of $j$-ary probability
distributions, thought of as the $(j - 1)$-simplex, $\Delta_{j - 1} \subset
\R^j$, and given the subspace topology. If $(p_i)_{i \in j} \in
\mathcal{C}(j)$, and for $i \in \{1, \ldots, j\}$, $(q_{i l})_{l \in k_i} \in
\mathcal{P}(k_i)$, we let
\[ \gamma ((p_i)_{i \in j} \otimes (q_{1 l})_{l \in k_0} \otimes \cdots
   \otimes (q_{j l})_{l \in k_{j - 1}}) = (p_i q_{i l})_{l \in k_i, i \in j}
   \in \mathcal{C}(k_0 + \cdots + k_{j - 1}) . \]
Basically, this says that, given a binary variable $X \in (x_i)_{i \in j}$ with
probability distribution $(p_i)_{i \in j}$, we refine the possible values of
$X$, splitting up each $x_i$. 

\smallskip

As a heuristic description of this procedure, imagine we are 
measuring physical systems and have suddenly discovered how to measure spin or
some other quantity that we were ignorant of before. Now there are more
distinguishable states that we can measure. We know the probability
distribution of these new states given an old state $x_i$: it is $(q_{i l})_{l
\in k_i}$, corresponding to new distinguishable states $(x_{i l})_{l \in
k_i}$. Now $X \in (x_{i l})_{l \in k_i, i \in j}$ may take any of $k_0 +
\cdots + k_{j - 1}$ values and must have the probability distribution $(p_i
q_{i l})_{l \in k_i, i \in j}$. We see the unit in this operad is the unique
probability distribution $(1) \in \mathcal{P}(1)$.

\smallskip
 
An important $\mathcal{P}$-algebra in ${\rm Cat} ({\rm Top})$ is
given by the additive monoid $\R_{\geqslant 0}$. As a category,
$\R_{\geqslant 0}$ is regarded as the one object category. the operad
$\mathcal{P}$ acts trivially on objects since there is only one object. On
maps, that is, on real numbers, we have 
\[ (p_i)_{i \in j} \cdot (x_i)_{i \in j} = \sum_i p_i x_i. \]

Since $\mathcal{P}$-algebras $A$ are also categories, we
can define an internal $\mathcal{P}$-algebra in $A$ as a lax map $1
\rightarrow A$ of $\mathcal{P}$-algbras where $1$ is the terminal
$\mathcal{P}$-algebra in ${\rm Cat}$ (see \cite{BaFriLei}, \cite{BaFriLei2})  
for details). This basically is an object $a \in A$ and, for each 
$p \in \mathcal{P}(j)$, a map $\alpha_p : p (a, \ldots, a) \rightarrow a$ such that
\[ \alpha_{p \circ (q_1, \ldots, q_n)} = \alpha_p \circ p (\alpha_{q_1},
   \ldots, \alpha_{q_n})   \ \ \   \text{ for every }  p \in \mathcal{P}(n) \ \text{ and } \
   q_i \in \mathcal{P}(m_i) \]
\[ \alpha_{\sigma p} = \alpha_p \ \ \  \text{ for every } \ \ 
p \in \mathcal{P}(n) \ \text{ and } \  \sigma \in {\rm Sym}_n \]
\[ \alpha_1 = \eta . \]

For $\R_{\geqslant 0}$, there is only one object, so $a
=\R_{\geqslant 0}$, and $\alpha$ is a map taking probability
distributions to positive real numbers satisfying the following four axioms:

\begin{enumerate}
\item For every $p \in \mathcal{P}(n)$ and $q_i \in \mathcal{P}(m_i)$
\[ \alpha (p \circ (q_1, \ldots, q_n)) = \alpha (p) + \sum_i p_i \alpha (q_i); \]
\item $\alpha ((1)) = 0$;
\item for every $p \in \mathcal{P}(n)$ and $\sigma \in {\rm Sym}_n$
\[  \alpha (\sigma p) = \alpha (p) \]
\item $\alpha : \mathcal{P}(n) \rightarrow \R_{\geqslant 0}$ is
  continuous for all $n$.
\end{enumerate}

\smallskip

Note that, in the first of these, composition of maps in the one object
category $\R_{\geqslant 0}$ is addition of real numbers. We require
the last one since we are looking for functoriality in ${\rm Cat}
({\rm Top})$. As it turns out (see \cite{BaFriLei}, \cite{BaFriLei2}), by 
Fadeev's theorem \cite{Fad}, the only function
satisfying these axioms, up to positive scalar multiples, is the Shannon
entropy, ${\rm Sh}$.

\smallskip

\subsection{Binary guessing trees}
Consider now a general binary information measure, $S : [0, 1] \rightarrow
\R_{\geqslant 0}$. We will assume that $S$ satisfies the identity
axioms, so that we can keep our approach finite rather than full of infinite
amounts of trivial flotsam. We can build an information measure on ternary
variables in several ways. For example, if we are trying to guess at the value
of $X$, which we know must be in $\{x_1, x_2, x_3 \}$, using only yes-or-no
questions, we could employ one of the following two strategies:

\begin{enumerate}
\item Is $X = x_1$? If not, is $X = x_2$?
\item Is $X = x_1$ or $x_2$? If yes, is $X = x_1$?
\end{enumerate}

Indeed, we see that any strategy that avoids asking trivial or irrelevant
questions arises as one of these strategies with a permutation of $\{1, 2,
3\}$. This gives us $2 \cdot 3! = 12$ possible ternary information measures.
There is a useful way of parametrizing these guessing strategies with rooted
trees.

\begin{prop}\label{bintreeprop}
  Let $S$ be a binary information measure with identity. For each $n \geqslant
  2$, there is a one-to-one correspondence between rooted full binary trees
  with $n$ leaves with labels in $\{1, \ldots, n\}$ and $n$-ary information
  measures arising from $S$.
\end{prop}

\proof Let ${\bf T}$ be a tree as above. We call such a tree an $(n,
2)$-tree. What it means to be full is that every vertex is either a leaf or
has two children. We will see that eliminating the single-child nodes is
equivalent to eliminating trivial and irrelevant questions from our set of
possible questions, making it finite. To see what is the set of possible questions,
consider that, if at a certain time we are certain that $X \in \{x_1, \ldots, x_n \}$, 
the yes-or-no questions available to us are exactly those of the form ``is $X \in
A$?'', where $A$ is a subset of $\{x_1, \ldots, x_n \}$. We label the leaves
of ${\bf T}$ with the possible values of $X$ according to their original
labels ($i \mapsto x_i$). The vertices which are not leaves are uniquely
labeled with the list of $x_i$ which label the leaves of their subtree. The
vertices will represent states of our knowledge of $X$ in that the labels will
denote the possible values of $X$ given what we have already measured.
Naturally, we begin at the root vertex, sure only that $X$ is one of the
$x_i$. At any vertex which is not a leaf, there are two child subtrees: a left
one, ${\bf L}$, and a right one, ${\bf D}$. Let $L$ be the set of
leaf labels of ${\bf L}$, $D$ those of ${\bf D}$. The true value of
$X$ must lie in either $L$ or $D$. Our question then is ``is $X \in L$?''. If
the answer is yes, we move to the left child. If the answer is no, we move to
the right child. At a leaf, we have ruled out all the possible values of $X$
except the one labeling our current vertex. 

As an example, consider the rooted full binary tree in Figure \ref{figtree1}:

\begin{figure}[h]
  \includegraphics{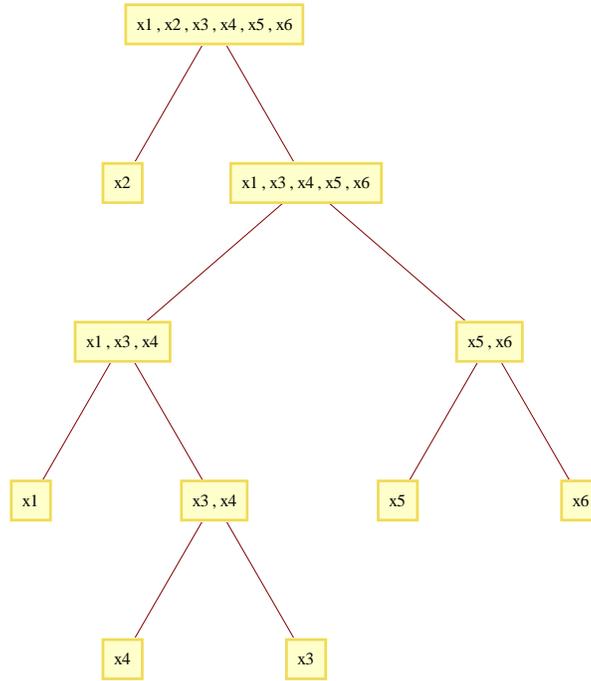}
  \caption{A rooted full binary tree \label{figtree1}} 
\end{figure}

We see $X$ lies in $\{x_1, \ldots, x_6 \}$. Our strategy goes like this:

1. Our first question is ``is $X = x_2$?''.

\ \ \ \ 1.1. If yes, we are done; $X = x_2$.

\ \ \ \ 1.2. If no, we ask ``is $X \in \{x_1, x_3, x_4 \}$?''.

\ \ \ \ \ \ \ \ \ 1.2.1. If yes, we ask ``is $X = x_1$?''.

\ \ \ \ \ \ \ \ \ \ \ \ \ \ 1.2.1.1. If yes, we are done; $X = x_1$.

\ \ \ \ \ \ \ \ \ \ \ \ \ \ 1.2.1.2. If no, we ask ``is $X = x_4$?''.

\ \ \ \ \ \ \ \ \ \ \ \ \ \ \ \ \ \ \ 1.2.1.2.1. If yes, we are done; $X =
x_4$.

\ \ \ \ \ \ \ \ \ \ \ \ \ \ \ \ \ \ \ 1.2.1.2.2. If no, we are also done; $X
= x_3$.

\ \ \ \ \ \ \ \ \ 1.2.2. If no, we ask ``is $X = x_5$?''.

\ \ \ \ \ \ \ \ \ \ \ \ \ \ 1.2.2.1. If yes, we are done; $X = x_5$.

\ \ \ \ \ \ \ \ \ \ \ \ \ \ 1.2.2.2. If no, we are also done; $X = x_6$.

\smallskip

Suppose these possible values occur with probabilities $p_1, \ldots, p_6$. We
see that the information measure corresponding to the above tree is
\[ S (p_2) + (1 - p_2) S ( \frac{p_1 + p_4 + p_3}{1 - p_2}) + (p_1 + p_4 + p_3)
S ( \frac{p_1}{p_1 + p_4 + p_3})  \]
\[ + (p_4 + p_3) S ( \frac{p_4}{p_4 + p_3}) +
(p_5 + p_6) S ( \frac{p_5}{p_5 + p_6}). \]
Note that permuting the labels of the leaves permutes the $p_i$.

\smallskip

Conversely, since any question is of the form ``is $X \in A$?'' for some
subsets $A$, we can build our tree inductively identifying $A$ with $L$ at a
given vertex, and labeling with the possible values of $X$ as we go,
beginning with the root. Any guessing strategy must exhaust the possibilities
for $X$, so any tree constructed in this way will be a well-defined $(n,
2)$-tree. Clearly this is the inverse process to the one described above. As
an example, suppose we want to guess at an $X \in \{x_1, x_2, x_3, x_4, x_5
\}$. First we might ask if $X \in \{x_1, x_2, x_4 \}$. If yes, we could ask if
$X = x_1$. If not, if $X = x_2$. Backtracking, if $X \notin \{x_1, x_2, x_4 \}$,
we could ask whether $X = x_5$. This strategy exhausts the possibilities for
$X$. It is represented by the tree in Figure \ref{figtree2}:

\begin{figure}[h]
  \includegraphics{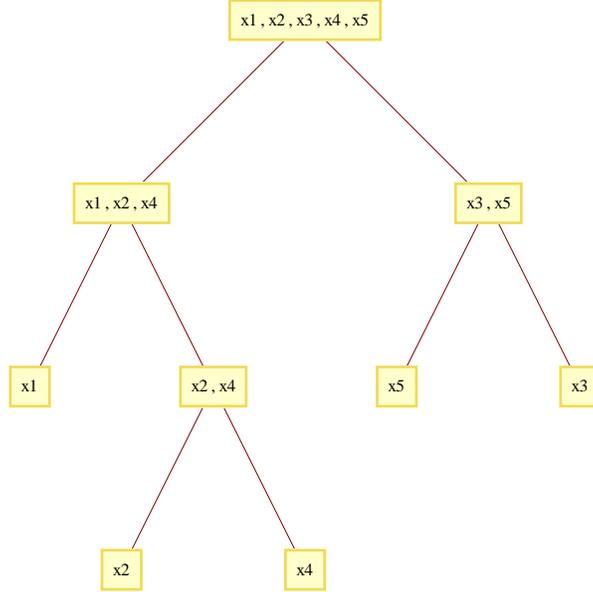}
    \caption{A guessing strategy \label{figtree2}}
\end{figure}

\smallskip

Given an $(n, 2)$-tree ${\bf T}$, there is a canonical way of
arranging and parenthesizing an expression of the form $x_1 \oplus_S \cdots
\oplus_S x_n$ so that it may be evaluated. This is the same one given in the
Catalan number identity \cite{Dav}. We consider the tree ${\bf T}'$ which is
labelled $1, \ldots, n$ from left to right. Let $\sigma_{{\bf T}} \in
{\rm Sym}_n$ be the permutation that sends the left-to-right labeling to the
original one on ${\bf T}$. We define $(x_1 \oplus_S \cdots \oplus_S
x_n)_{{\bf T}} = (x_{\sigma_{{\bf T}} (1)} \oplus_S \cdots \oplus_S
x_{\sigma_{{\bf T}} (n)})_{{\bf T}'}$. Thus, it suffices to consider
the case when ${\bf T}$ is labelled left-to-right. In this case, there is
a $1 \leqslant r < n$ such that for $1 \leqslant j \leqslant r$, $x_j$ is a
label of a leaf of the left subtree ${\bf L}$, ie. $x_j \in L$, and for
all $r < j < n$, $x_j \in D$, where ${\bf D}$ is the right subtree of
${\bf T}$. Then we define inductively $(x_1 \oplus_S \cdots \oplus_S
x_n)_{{\bf T}} = (x_1 \oplus_S \cdots \oplus_S x_r)_{{\bf L}} + (x_{r +
1} \oplus_S \cdots \oplus_S x_n)_{{\bf D}}$, with a tree with two children
${\bf T}_2$ giving $(x_1 \oplus_S x_2)_{{\bf T}_2} = x_1 \oplus_S x_2$.

\begin{thm}
  Given an $(n, 2)$-tree, ${\bf T}$, and a binary information measure $S$
  with identity, the following holds:
  \[ (x_1 \oplus_S \cdots \oplus_S x_n)_{{\bf T}} = \min_{\sum p_i = 1} ( \sum
  p_i x_i - T S_{{\bf T}} (p_1, \ldots, p_n)). \]
\end{thm}

\proof Before we begin the proof in ernest, we illustrate the argument with an explicit 
example. We
see the tree in Figure \ref{figtree1}, which we denote ${\bf T}$, corresponds to the
arrangement of parentheses $x_1 \oplus_S ((x_2 \oplus_S (x_3 \oplus_S x_4)) \oplus_S
(x_5 \oplus_S x_6))$ and the permutation $\sigma = (1 2) (3 4) \in
{\rm Sym}_6$. We calculate

\[  x_1 \oplus_S ((x_2 \oplus_S (x_3 \oplus_S x_4)) \oplus_S (x_5 \oplus_S x_6)) = \]
\[ \min_{p_1} (p_1 x_1 + (1 - p_1) ((x_2 \oplus_S (x_3 \oplus_S x_4)) \oplus_S (x_5
\oplus_S x_6)) - T S (p_1)) \]
\[ = \min_{p_1} (p_1 x_1 + (1 - p_1) \min_{p_2} (p_2 (x_2 \oplus_S (x_3 \oplus_S
x_4)) + (1 - p_2) (x_5 \oplus_S x_6) - T S (p_2)) - T S (p_1)) \]
\[ = \min_{p_1, p_2} \left( p_1 x_1 
+ (1 - p_1) p_2  \min_{p_3} \left( p_3 x_2 + (1 - p_3)
(x_3 \oplus_S x_4) - T S (p_3)\right)  \right. \]
\[ \left. + (1 - p_1) (1 - p_2) \min_{p_4} \left( p_4 x_5
+ (1 - p_4) x_6 - T S (p_4)\right) 
- T (S (p_1) + (1 - p_1) S (p_2)) \right) \]
\[ = \min_{p_1, p_2, p_3, p_4, p_5} ( p_1 x_1 + (1 - p_1) p_2 p_3 x_2 + (1 - p_1)
p_2 (1 - p_3) p_5 x_3 \] 
\[ + (1 - p_1) p_2 (1 - p_3) (1 - p_5) x_4
+ (1 - p_1) (1 - p_2) p_4 x_5 + (1 - p_1) (1 - p_2) (1 - p_4) x_6 \]
\[ - T ( S (p_1) 
+ (1 - p_1) S (p_2) + (1 - p_1) p_2 S (p_3) 
+ (1 - p_1) (1 - p_2) S (p_4) + (1 - p_1) p_2 (1 - p_3) S (p_5)   ))    . \]

Now we make the substitution
\begin{align*}
q_1 = & p_1 \\
q_2 = & (1 - p_1) p_2 p_3 \\
q_3 = & (1 - p_1) p_2 (1 - p_3) p_5 \\
q_4 = & (1 - p_1) p_2 (1 - p_3) (1 - p_5) \\
q_5 = & (1 - p_1) (1 - p_2) p_4 \\
q_6 = & (1 - p_1) (1 - p_2) (1 - p_4) .
\end{align*}

We notice $q_1 + \cdots + q_6 = 1$, and
\begin{align*}
p_1 = &  q_1 \\
p_2 = & (q_2 + q_3 + q_4) / (1 - q_1) \\
p_3 = & q_2 / (q_2 + q_3 + q_4) \\
p_4 = & q_5 / (q_5 + q_6) \\
p_5 = & q_3 / (q_3 + q_4) . 
\end{align*}

Notice that these look like relative probabilities. This is no coincidence.
Making this substitution above yields
\[ x_1 \oplus_S ((x_2 \oplus_S (x_3 \oplus_S x_4)) \oplus_S (x_5 \oplus_S x_6)) = \]
\[  \min_{\sum q_i = 1} ( \sum q_i x_i
- T (S (q_1) + (1 - q_1) S ( \frac{q_2 + q_3 + q_4}{1 - q_1}) \] 
\[ + (q_2 + q_3 +
q_4) S ( \frac{q_2}{q_2 + q_3 + q_4}) + (q_3 + q_4) S ( \frac{q_3}{q_3 + q_4})
+ (q_5 + q_6) S ( \frac{q_5}{q_5 + q_6})) . \]

Applying $\sigma$ we obtain 
\[ (x_1 \oplus_S x_2 \oplus_S x_3 \oplus_S x_4 \oplus_S x_5
\oplus_S x_6)_{{\bf T}} = \min_{\sum p_i = 1} ( \sum p_i x_i - T
S_{{\bf T}} (p_1, \ldots, p_6)) \] 
as the theorem claims. 

Now we are ready to prove the theorem in general.

\begin{lem}\label{treelem1}
  Suppose that, at the root, the tree ${\bf T}$ has left subtree ${\bf L}$ with $l$
  leaves, and right subtree ${\bf D}$ with $d$ leaves, and the leaves of
  ${\bf T}$ are labeled left to right. Then
  \[  S_{{\bf T}} (p_1, \ldots, p_l, p_{l + 1}, \ldots, p_{l + d}) = \]
  \[  S (p_1 + \cdots + p_l)
   + (p_1 + \cdots + p_l) S_{{\bf L}} ( \frac{p_1}{p_1 + \cdots + p_l},
  \ldots, \frac{p_l}{p_1 + \cdots + p_l}) \]
  \[ + (p_{l + 1} + \cdots + p_{l + d}) S_{{\bf D}} ( \frac{p_{l + 1}}{p_{l
  + 1} + \cdots + p_{l + d}}, \ldots, \frac{p_{l + d}}{p_{l + 1} + \cdots +
  p_{l + d}}) . \]
\end{lem}

\begin{lem}\label{treelem2}
  Suppose at the root, ${\bf T}$ has left subtree ${\bf L}$ with $l$
  leaves, and right subtree ${\bf D}$ with $d$ leaves, and the leaves of
  ${\bf T}$ are labeled left to right. Then
  \[  (x_1 \oplus_S \cdots \oplus_S x_l \oplus_S x_{l + 1} \oplus_S \cdots \oplus_S x_{l +
  d})_{{\bf T}} \]
  \[  = \min_{p} (p (x_1 \oplus_S \cdots \oplus_S x_l)_{{\bf L}} + (1 - p)
  (x_{l + 1} \oplus_S \cdots \oplus_S x_{l + d})_{{\bf D}} - T S (p)) .\] 
\end{lem}

The proof of both of these statements is immediate from the definitions.

\smallskip

Now, clearly the theorem holds when ${\bf T}$ has two leaves, and
since our trees are full, we can use this as the base case in an induction.

\smallskip

Suppose the theorem holds for all trees with less than $n$ leaves. Let
${\bf T}$ be an $(n, 2)$-tree with leaves labeled from left to right. At
the root, since ${\bf T}$ is full, ${\bf T}$ has nonempty left and
right subtrees, ${\bf L}$ and ${\bf D}$, with $l > 0$ and $d > 0$
leaves, respectively, such that $l + d = n$, so $l, d < n$. By the inductive
hypothesis and the second lemma above,
\[  (x_1 \oplus_S \cdots \oplus_S x_n)_{{\bf T}} = \min_p (p \min_{p_1 + \cdots
+ p_l = 1} ( \sum p_i x_i - T S_{{\bf L}} (p_1, \ldots, p_l)) \]
\[  + (1 - p) \min_{p_{l + 1} + \cdots + p_{l + d} = 1} ( \sum p_i x_i - T
S_{{\bf D}} (p_{l + 1}, \ldots, p_{l + d})) - T S (p)) . \]

Make the substitution $q_i = p p_i$, for each $i \in \{1, \ldots, l\}$, and $q_i
= (1 - p) p_i$, for each $i \in \{l + 1, \ldots, l + d\}$. Note that $q_1 + \cdots +
q_l = p$ and $q_{l + 1} + \cdots + q_{l + d} = 1 - p$. This yields
\[  (x_1 \oplus_S \cdots \oplus_S x_n)_{{\bf T}} = \min_{\sum q_i = 1} ( \sum
q_i x_i \]
\[ - T ((q_1 + \cdots + q_l) S_{{\bf L}} ( \frac{q_1}{q_1 + \cdots + q_l},
\ldots, \frac{q_l}{q_1 + \cdots + q_l}) \]
\[ + (q_{l + 1} + \cdots + q_{l + d}) S_{{\bf D}} ( \frac{q_{l + 1}}{q_{l +
1} + \cdots + q_{l + d}}_{}, \ldots, \frac{q_{l + d}}{q_{l + 1} + \cdots +
q_{l + d}}_{}) \]
\[ + S (q_1 + \cdots + q_l))), \]
which by the first lemma is
\[ \min_{\sum q_i = 1} ( \sum q_i x_i - T S_{{\bf T}} (q_1, \ldots, q_n))). \]

\smallskip

We need now show that this holds for arbitrary labelings of the leaves of
${\bf T}$. If $\sigma$ is a permutation of $\{1, \ldots, n\}$, then
\[  (x_{\sigma (1)} \oplus_S \cdots \oplus_S x_{\sigma (n)})_{{\bf T}} =
\min_{\sum q_i = 1} ( \sum q_i x_{\sigma (i)} - T S_{{\bf T}} (q_1,
\ldots q_n)) \]
\[  = \min_{\sum p_i = 1} ( \sum p_i x_i - T S_{{\bf T}} (p_{\sigma (1)},
\ldots, p_{\sigma (n)})), \]
where we have substituted $p_i = q_{\sigma^{- 1} (i)}$. This proves the
theorem.
\endproof

The connection between these guessing games and the thermodynamics of mixing discussed in \ref{statmechSec} can be intuited in the following way. The entropy of a system arises from considering the ``correct counting" of states. In more words, some states are indistinguishable from others, and this affects their multiplicity in the partition sum. The entropy function tells us what the overall degree of distinguishability is. We can see this point of view in Boltzmann's famous equation asserting $S = k_B \log \Omega$, where $\Omega$ is the number of microstates which degenerate to a given macrostate. When we perform mixtures in a certain order, we are giving an order to this distinguishing process, as we are when we decide on an order to ask questions in a guessing game, distinguishing possible values from impossible values of the unknown variable.

\subsection{General guessing trees}

Now suppose for each $n \in V \subseteq \{m \in \N\,|\,m
\geqslant 2\}$ we have an $n$-ary information measure $S_n$. We want to impose
the following condition.

\begin{enumerate}

\item  (Coherence) Suppose $n > m$ and, for all but $1 < i_1 < \cdots < i_m < n$,
  $p_j = 0$. Then
  \[   S_n (p_1, \ldots, p_n) = S_m (p_{i_1}, \ldots, p_{i_m}).   \]
\end{enumerate}

We can always write 
\[ S_{n - 1} (p_1, \ldots, p_{n - 1}) = S_n (p_1, \ldots,
p_{n - 1}, 0), \] 
so that we can take $V$ as an initial segment of
$\N_{\geqslant 2}$. This way, we can instead think about $v = \sup
V$. Many definitions of entropies have $v = \infty$. Examples include the
Shannon, Renyi, and Tsallis entropies. These are generally defined by
functions $f, g$ such that
\[ S_n (p_1, \ldots, p_n) = f ( \sum_{1 \leqslant i \leqslant n} g (p_i)) . \]
Any entropy of this form trivially satisfies the coherence axiom.

\smallskip

In this more general setting, we can ask any question with up to $v$
possible answers. This potentially gives us many new ways to play guessing
games, or equivalently, to build more general information measures. For
example, if $v \geqslant 5$, then we can measure a 12-ary random variable $X
\in \{x_1, \ldots, x_{12} \}$ by asking first whether $X \in \{x_1, \ldots,
x_5 \}$. If yes, we simply measure the value of $X$. Otherwise, ask which of
$\{x_6, x_7, x_8 \}$, $\{x_9, x_{12} \}$, or $\{x_{10}, x_{11} \}$ contains
$X$, and then simply measure the value of $X$ (note that order may matter:
$S_n$ may not be symmetric). This gives us information
\[  S_2 (p_1 + \cdots + p_5, p_6 + \cdots + p_{12})  \] 
\[ + (p_1 + \cdots + p_5) S_5 (
\frac{p_1}{p_1 + \cdots + p_5}, \ldots, \frac{p_5}{p_1 + \cdots + p_5})  \]
\[ + (p_6 + \cdots + p_{12}) S_3 ( \frac{p_6 + p_7 + p_8}{p_6 + \cdots + p_{12}},
\frac{p_9 + p_{12}}{p_6 + \cdots + p_{12}}, \frac{p_{10} + p_{11}}{p_6 +
\cdots + p_{12}})  \]
\[ + (p_6 + p_7 + p_8) S_3 ( \frac{p_6}{p_6 + p_7 + p_8}, \frac{p_7}{p_6 + p_7 +
p_8}, \frac{p_8}{p_6 + p_7 + p_8}) \]
\[ + (p_9 + p_{12}) S_2 ( \frac{p_9}{p_9 + p_{12}}, \frac{p_{12}}{p_9 + p_{12}}) \]
\[ + (p_{10} + p_{11}) S_2 ( \frac{p_{10}}{p_{10} + p_{11}}, \frac{p_{11}}{p_{10}
+ p_{11}}), \]
where now we write $S_2$ as a two-variable function for
consistency. We see something extremely similar to the binary 
case is happening here.

\begin{prop}
  Let $n, v \geqslant 2$, and suppose for each $2 \leqslant j < v + 1$ we have
  a $j$-ary information measure $S_j$, and these together satisfy the
  coherence axiom. Guessing strategies of $n$-ary random variables where we
  allow questions of up to $v$ possible answers are in bijective
  correspondence with the set of $(n, v)$-trees, rooted trees with labelled 
  leaves such that every vertex is either a leaf or has between $2$ and $v$
  children.
\end{prop}

\proof Every relevant question that can be asked is of the form ``which of
$A_1, \ldots, A_m$ contains $X$?'' for certain disjoint subsets $A_1, \ldots,
A_m$. We identify these subsets with the leaves of the $m$ subtrees extending
from the current vertex, once again identifying the vertices with states of
our knowledge of $X$. 

\smallskip

\begin{figure}[h]
 \includegraphics{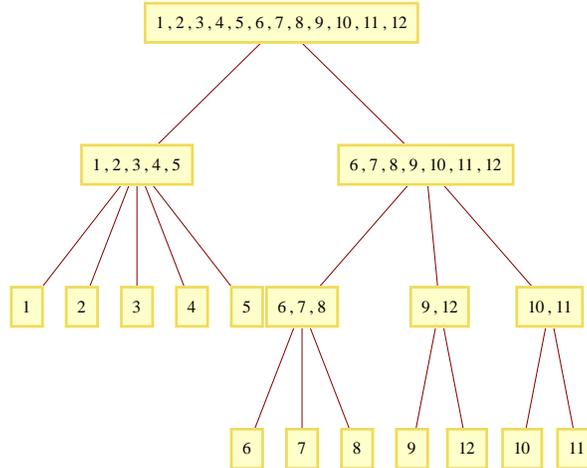}
 \caption{An $(n, v)$-tree \label{figtree3}}
\end{figure}

For example, from the previous algorithm we have the tree in Figure \ref{figtree3}.

Conversely, to go from an $(n, v)$-tree to a guessing strategy one must only
follow the tree to its leaves.

\smallskip

Now we have some basic $n$-ary functions for more than just $n = 2$.
Namely, we can define
\begin{equation}\label{naryfunction}
  x_1 \oplus_S \cdots \oplus_S x_n := \min_{\sum p_i = 1} ( \sum_{1
  \leqslant i \leqslant n} p_i x_i - T S_n (p_1, \ldots, p_n)) .
\end{equation}

This has a thermodynamic interpretation of a simultaneous mixing of $n$ gas species. We have the following result.

\begin{prop}
  Let $n > 2$. The following hold.
 \begin{enumerate}
  \item For every $j$ 
  \[ x_1 \oplus_S \cdots \oplus_S x_j \oplus_S \infty \oplus_S x_{j + 2} \oplus_S
  \cdots \oplus_S x_n = \] \[  x_1 \oplus_S \cdots \oplus_S x_j \oplus_S x_{j + 2} \oplus_S
  \cdots \oplus_S x_n \] if and only if the $S_n$ share the coherence property.
  For $n = 2$ this is the identity property.
  
  \item $x_1 \oplus_S \cdots \oplus_S x_n$ is symmetric if and only if 
  $S_n$ is symmetric. For $n = 2$ this is commutativity.
\end{enumerate}  
\end{prop}

The proof of this fact is immediate.

\smallskip

We can generalize the parentheses correspondence to $(n, v)$-trees:
to any $(n, v)$-tree ${\bf T}$ we can associate a unique $n$-ary function
$(x_1 \oplus_S \cdots \oplus_S x_n)_{{\bf T}}$ given by arranging parentheses
according to ${\bf T}$ and $x_i$ according to the labels on the leaves.

\smallskip

For example, to the tree above we associate
\[ (x_1 \oplus_S x_2 \oplus_S x_3 \oplus_S x_4 \oplus_S x_5) \oplus_S ((x_6
   \oplus_S x_7 \oplus_S x_8) \oplus_S (x_9 \oplus_S x_{12}) \oplus_S (x_{10} \oplus_S
   x_{11})) . \]
   
We have then the analog for $(n,v)$-trees of Lemmas \ref{treelem1} 
and \ref{treelem2}.
   
\begin{lem}\label{ntreelem1}
  Suppose the root of an $(n, v)$-tree ${\bf T}$ has sub-$(l_j, v)$-trees
  (resp. from left to right) ${\bf A}_1, \ldots, {\bf A}_m$, and the
  leaves of ${\bf T}$ are labeled left to right. Let $L_j = l_1 + \cdots
  + l_j$ and $L_0 = 0$. Then the following holds.
  \[ S_{{\bf T}} (p_1, \ldots, p_n) = \] \[ \sum_{1 \leqslant j \leqslant m}
     (p_{L_{j - 1} + 1} + \cdots + p_{L_j}) S_{{\bf A}_j} ( \frac{p_{L_{j
     - 1} + 1}}{p_{L_{j - 1} + 1} + \cdots + p_{L_j}}, \ldots,
     \frac{p_{L_j}}{p_{L_{j - 1} + 1} + \cdots + p_{L_j}}) \]
  \[ + S_m (p_{L_0 + 1} + \cdots + p_{L_1}, \ldots, p_{L_{m - 1} + 1} + \cdots
     + p_{L_m}) \]
\end{lem}

\medskip

\begin{lem}\label{ntreelem2}
  Suppose the root of an $(n, v)$-tree ${\bf T}$ has sub-$(l_j, v)$-trees
  (resp. from left to right) ${\bf A}_1, \ldots, {\bf A}_m$, and the
  leaves of ${\bf T}$ are labeled left to right. Let $L_j = l_1 + \cdots
  + l_j, L_0 = 0$. Then the following holds:
  \begin{align*}  (x_1 \oplus_S \cdots \oplus_S x_n)_{{\bf T}} =& \min_{\sum q_i =
     1} (q_1 (x_1 \oplus_S \cdots \oplus_S x_{l_1})  \\ + & \cdots + q_m (x_{l_1 +
     \cdots + l_{m - 1} + 1} \oplus_S \cdots \oplus_S x_{l_1 + \cdots + l_m}) \\ - & T
     S_m (q_1, \ldots, q_m)) . 
   \end{align*}  
\end{lem}

As before, both of these are immediate from the definitions. Finally, we have
the theorem:

\begin{thm}\label{ntreethm}
  Given an $(n, v)$-tree ${\bf T}$, and for each $2 \leqslant j
  \leqslant n$ an information measure $S_j$, such that together they satisfy
  the coherence axioms, the following holds:
  \[  (x_1 \oplus_S \cdots \oplus_S x_n)_{{\bf T}} = \min_{\sum p_i =
     1} ( \sum p_i x_i - T S_{{\bf T}} (p_1, \ldots, p_n)) . \]
\end{thm}

\proof Once again we proceed by strong induction on the number of leaves. We
know the theorem holds for $n = 2$. Suppose the theorem holds for every $(m,
v)$-tree with $m < n$. Let ${\bf T}$ be an $(n, v)$-tree with leaves
labeled from left to right. ${\bf T}$ has $k \geqslant 2$ sub-$(l_i,
v)$-trees starting at the root (resp. from left to right) ${\bf A}_1,
\ldots, {\bf A}_k$ with $l_i > 0$. We must have $l_1 + \cdots + l_k = n$,
so $l_i < n$. By the inductive hypothesis and the second lemma above then,
\[ (x_1 \oplus_S \cdots \oplus_S x_n)_{{\bf T}} = \min_{\sum q_i = 1} (q_1
\min_{p_1 + \cdots + p_{l_1} = 1} (p_1 x_1 + \cdots + p_{l_1} x_{l_1} - T
S_{{\bf A}_1} (p_1, \ldots, p_{l_1}) + \cdots \]
\[  + q_k \min_{p_{l_1 + \cdots + l_{k - 1} + 1} + \cdots + p_{l_1 + \cdots + l_k}
= 1} ( \sum_{j = l_1 + \cdots + l_{k - 1} + 1}^{l_1 + \cdots + l_k} p_j x_j -
T S_{{\bf A}_k} (p_{l_1 + \cdots + l_{k - 1} + 1}, \ldots, p_{l_1 +
\cdots + l_k})) \]
\[ - T S_k (p_1 + \cdots + p_{l_1}, \ldots, p_{l_1 + \cdots + l_{k - 1} + 1} +
\cdots + p_n)) . \]

\smallskip

For each $i \in \{1, \ldots, k\}$, and each $j \in \{l_1 + \cdots + l_{i - 1}
+ 1, \ldots, l_1 + \cdots + l_i \}$, where we define $l_0 = 0$, we make the
substitution $\tilde{q}_j = q_i p_j$. That way, $\sum_{j = l_1 + \cdots + l_{i
- 1} +}^{l_1 + \cdots + l_i} \tilde{q}_j = q_i,$ so we have
\[  (x_1 \oplus_S \cdots \oplus_S x_n)_{{\bf T}} = \min_{\sum \tilde{q}_j = 1} (
\sum \tilde{q}_j x_j - T (( \tilde{q}_1 + \cdots + \tilde{q}_{l_1})
S_{{\bf A}_1} ( \frac{\tilde{q}_1}{\tilde{q}_1 + \cdots +
\tilde{q}_{l_1}}, \ldots) + \cdots   \]
\[ + ( \tilde{q}_{l_1 + \cdots + l_{k - 1} + 1} + \cdots + \tilde{q}_n)
S_{{\bf A}_k} ( \frac{\tilde{q}_{l_1 + \cdots + l_{k - 1} +
1}}{\tilde{q}_{l_1 + \cdots + l_{k - 1} + 1} + \cdots + \tilde{q}_n}, \ldots)) . \]

By Lemma \ref{ntreelem1}, this equals
\[ \min_{\sum \tilde{q}_j = 1} ( \sum \tilde{q}_j x_j - T S_{{\bf T}} (
\tilde{q}_1, \ldots, \tilde{q}_n)).  \]

Now let $\sigma$ be any permutation of $\{1, \ldots, n\}$. We see

\[ (x_{\sigma (1)} \oplus_S \cdots \oplus_S x_{\sigma (n)})_{{\bf T}} =
\min_{\sum q_i = 1} ( \sum q_i x_{\sigma (i)} - T S_{{\bf T}} (q_1,
\ldots q_n)) \]
\[ = \min_{\sum p_i = 1} ( \sum p_i x_i - T S_{{\bf T}} (p_{\sigma (1)},
\ldots, p_{\sigma (n)})) , \]
where we have substituted $p_i = q_{\sigma^{- 1} (i)}$. This proves the
theorem.
\endproof

\subsection{Information algebra}\label{infoalg}

We define $\mathcal{T}_v (n)$ to be the class of $(n, v)$-trees such that
$\mathcal{T}_v (0)$ contains only the empty graph and $\mathcal{T}_v (1)$
contains only the unique one-leav\`ed $(n, v)$-tree. We put an operad structure on the union of these collections, $\mathcal{T}$ with
composition given by leaf-to-root composition of trees, which is clearly unital,
associative, and ${\rm Sym}$-equivariant. Our underlying category is the cartesian monoidal category of sets of graphs, with $\kappa=\mathcal{T}_v (1)$. Note that this unital operad structure, if also given a free group structure, forms the well-known $A_\infty$-operad.

\smallskip

Consider the one-object topological category, $R$, and a coherent set
\[ \{S_j : I^j \rightarrow \R_{\geqslant 0} \, | \, 2 \leqslant j < v + 1\} \]
of information measures. For each $n \geqslant 2$, and each ${\bf T} \in
\mathcal{T} (n)$, we define
\[ {\bf T}(x_1, \ldots, x_n) = \min_{\sum p_i = 1} ( \sum p_i x_i - T
   S_{{\bf T}} (p_1, \ldots, p_n)) . \]
With the definition of $S_{{\bf T}}$, as in the previous section. For
$\mathcal{T} (1)$, we define $(x)_{{\bf T}} = x$, for $\mathcal{T}_v
(0) = \kappa$, we define $()_{{\bf T}} = \infty$.

\smallskip

By Theorem \ref{ntreethm} above,
this is the same as $(x_1 \oplus_S \cdots \oplus_S x_n)_{{\bf T}}$, which
clearly behaves well under composition of trees and the action of
${\rm Sym}_n$, so this makes $R$ a $\mathcal{T}$-algebra, which we call the
{\em information algebra} $(R, S)$. This characterizes the complete algebraic structure of the Witt semiring $R$ over $K$ arising from $S$. The next proposition is written in the original
convention for semifields and summarizes some characteristics of this action,
each of which are immediate from the definitions.

\begin{prop}\label{infoalgstr}
  Let ${\bf T} \in \mathcal{T} (n)$, $x_1, \ldots, x_n, y \in K, \alpha
  \in \R_{\geqslant 0}$. Then the following hold.
 \begin{enumerate} 
  \item The $\mathcal{T}$-algebra structure on $R$ is additive: for all $1 \leqslant j \leqslant n$
  \[   {\bf T}(x_1, \ldots, x_{j - 1},
  x_j + y, x_{j + 1}, \ldots, x_n) = \] \[ {\bf T}(x_1, \ldots, x_j, \ldots,
  x_n) +{\bf T}(x_1, \ldots, y, \ldots, x_n).  \]
  
  \item Multiplication distributes over the $\mathcal{T}$-algebra structure:
  \[  y (x_1 \oplus_S \cdots \oplus_S x_n)_{{\bf T}} = (y x_1 \oplus_S \cdots
  \oplus_S y x_n)_{{\bf T}}. \]
  
  \item The $\mathcal{T}$-algebra structure also satifies  
  \[ (x_1 \oplus_S \cdots \oplus_S x_n)^{\alpha}_{{\bf T}} (T) = (x_1^{\alpha}
  \oplus_S \cdots \oplus_S x_n^{\alpha})_{{\bf T}} (\alpha T). \]
 \end{enumerate} 
\end{prop}

The relations which are most natural to consider are of the form
\[ {\bf T}_1 (x_1, \ldots, x_n) ={\bf T}_2 (x_1, \ldots, x_n) \ \ \ 
   \forall x_i \in R , \]
where ${\bf T}_1$ and ${\bf T}_2$ are $(n, v)$-trees acting on the
information algebra $(R, S)$. The reason is that we can interpret this as an
equivalence of guessing strategies, so these are exactly
the kind of relations that would define an information measure. Note it is not just ${\bf T}_1$ and
${\bf T}_2$ which are affected by this relation. Because composition of
trees gives the composition of their actions on $R$, whenever some tree can be
written ${\bf A} \circ {\bf T}_1 \circ ({\bf A}_1, \ldots,
{\bf A}_n)$, this is equivalent to ${\bf A} \circ {\bf T}_2 \circ
({\bf A}_1, \ldots, {\bf A}_n)$. The equivalence classes of these
trees for some fixed set of relations ${\bf R}$, defines a quotient
operad $\mathcal{T}/{\bf R}$ which is the set of possible guessing
strategies up to equivalence under the information measure. The terminal
object in this construction is the one with exactly one $(m, v)$-tree for each
$m$, which is precisely the quotient operad arising from the Shannon entropy.

\smallskip

However, one quickly finds that these simple relations are inadequate for describing 
the full range of information measures. If we have an equivalence of trees, we can always prune 
corresponding leaves by inserting the identity, $\infty$ in the current notation, in the place 
of that variable. For binary information measures, one has the followinng fact.

\begin{prop}\label{trivtreerel}
 Suppose $S$ is a commutative binary information measure, and ${\bf T}_1, {\bf T}_2$ are $(n, 2)$-trees.
 Either ${\bf T}_1 = {\bf T}_2$ is implied by commutativity or implies associativity, hence forces $S$ to
 be the Shannon entropy.
\end{prop}
\proof
We proceed by induction on $n$. When $n=3$, one checks the above is true by simply checking each case. 
Suppose the theorem holds for all 
$m < n$. By pruning a leaf, we see that either the 
relation implies associativity or the pruned subtrees are equal. In the case of the latter, prune a different leaf 
and we see the theorem holds.
\endproof

Thus, one may wish to pass into a setting where we may consider linear combinations of trees, ie. 
we can put a free vector space structure on our original operad. This gives us an $A_{\infty}$ operad, with the action on the information algebra $R$ extending uniquely under the 
Frobenius action.

\smallskip

Let us now consider what an internal $\mathcal{T}$-algebra in $R$ is.
For each $n$, this is a continuous map $\alpha_n : \mathcal{T}(n) \rightarrow
R$ such that the following hold:
\begin{enumerate}
\item for all ${\bf T} \in \mathcal{T}(n)$, and ${\bf A}_1 \in
   \mathcal{T}(m_1), \ldots, {\bf A}_n \in \mathcal{T}(m_n)$, 
\[ \alpha_{m_1 + \cdots + m_n} {\bf T} \circ ({\bf A}_1, \ldots,
   {\bf A}_n) = \alpha_n ({\bf T}) +{\bf T}(\alpha_{m_1}
   ({\bf A}_1), \ldots, \alpha_{m_n} ({\bf A}_n)); \]
   
\item for all ${\bf T} \in \mathcal{T}(n)$ and $\sigma \in {\rm Sym}_n$,   
\[  \alpha_n
   (\sigma {\bf T}) = \alpha_n ({\bf T}); \]

\item  The following condition also holds:  
\begin{equation}\label{alpha1T}
 \alpha_1 (\mathcal{T}(1)) = 0.
\end{equation} 
\end{enumerate}

To simplify notation, we will suppress 
the subscripts on $\alpha$ and just consider $\alpha :
\mathcal{T} \rightarrow \R_{\geqslant 0}$. 
The second condition above just
says that $\alpha ({\bf T})$ does not depend on the labels of
${\bf T}$. 

For each $2 \leqslant n < v + 1$ we define $h_n \in
\R_{\geqslant 0}$ as the unique value $\alpha$ takes on the $(n,
v)$-trees with $n + 1$ vertices, that is, those corresponding to $S_n$. 

Every tree in $\mathcal{T}$ is built from these basic trees, and by the first condition
above, so is $\alpha ({\bf T})$. 

If at the root ${\bf T}$ has
subtrees ${\bf A}_1, \ldots, {\bf A}_n$ from left to right, then
\[ \alpha ({\bf T}) = h_n + \alpha ({\bf A}_1) \oplus_S \cdots \oplus_S
   \alpha ({\bf A}_n) . \]
   
For the tree in Figure \ref{figtree3}, this gives
\[ h_2 + h_5 \oplus_S (h_3 + h_3 \oplus_S h_2 \oplus_S h_2) . \]

We see $+$ goes down the tree, and $\oplus_S$ goes across the tree. 

It is easy to see $h_3 \geqslant h_3 \oplus_S h_2 \oplus_S h_2$, and $h_2 \geqslant h_5
\oplus_S (h_3 \oplus_S h_2 \oplus_S h_2)$, so this can be simplified to
\[ h_5 \oplus_S (h_3 \oplus_S h_2 \oplus_S h_2), \]
which we see can be obtained through a different recursion strategy. Instead
of picking off the subtrees at the root, we could pick off the basic subtrees
just above the leaves. This is just another way of writing ${\bf T}$ as a
composition of trees, and puts the recursion into the second term rather than
the first in \eqref{alpha1T}. 

Because of this recursion, every internal
$\mathcal{T}$-algebra of $\R_{\geqslant 0}$ is determined by the
sequence $(h_j)_{2 \leqslant j < v + 1}$ (by the third condition above, 
implicity $h_1 = 0$). 

When $R =\R^{\max, \ast}_{\geqslant 0},$ and we use $S =
{\rm Sh}$, the Shannon entropy, then $x \oplus_S y = (x^{1 / T} + y^{1 / T})^T$, so 
the above becomes
\[ \alpha ({\bf T}) = \max (h_n, (\alpha ({\bf A}_1)^{1 / T} +
   \cdots + \alpha ({\bf A}_n)^{1 / T})^T) . \]

\section{Further perspectives and directions}

We sketch here some possible further directions where the notion of
thermodynamic semirings may prove useful.

\subsection{Information geometry}\label{infogeomsec}

Information geometry was developed \cite{Ama}, \cite{AmaNa}, \cite{IkTaAma}
as a way to encode, using methods based on Riemannian geometry,
statistical information, such as how to infer unobserved variables 
on the basis of observed ones by reducing conditional joint probabilities
to marginal distributions.

\smallskip

We consider a smooth univariate binary statistical $n$-manifold $\cQ$
as in Definition \ref{infomanifold} parameterized by $\eta \in \cX \subset \R^n$. One may deal with the multivariate case similarly. 
\smallskip

The Fisher information metric (see \cite{AmaNa})
on information manifolds is given by
\[ g_{i j} (\theta) = \int \frac{\partial \ln p (x ; \theta)}{\partial \theta_i}
\frac{\partial \ln p (x ; \theta)}{\partial \theta_j} d x \] 
and it defines a Riemannian metric on a statistical manifold $\cQ$.

\smallskip

Another important notion in information geometry is that of e-flat and m-flat submanifolds,
which we recall here.

\smallskip

A submanifold $\cS\subset \cQ$ is e-flat if, for all $t\in [0,1]$ and all $p(\eta)$ and
$q(\eta)$ in $\cS$ the mixture
$\log r(\eta,t) = t \log p(\eta) + (1-t) \log q(\eta) + c(t)$,
with $c(t)$ a normalization factor, is also in $\cS$.

\smallskip

A submanifold $\cS\subset \cQ$ is m-flat if, for all $t\in [0,1]$ and all $p(\eta)$ and
$q(\eta)$ in $\cS$ the mixture
$r(\eta,t)= t p(\eta)+ (1-t) q(\eta)$  
is also in $\cS$.

\smallskip

One-dimensional e-flat or m-flat manifolds are called e-geodesics and m-geodesics.
In information geometry, maximum posterior marginal optimization is achieved by
finding the point on an e-flat submanifold $\cS$ that minimizes the KL divergence, 
see \cite{AmaNa}, \cite{IkTaAma}. It turns out that the point on an e-flat submanifold 
$\cS$ that minimizes the KL divergence also minimizes the Riemannian metric 
given by the Fisher information metric.

\smallskip

More precisely, when considering the KL
divergences ${\rm KL}(p;q(\eta))$, where $q(\eta)$ varies in an
e-flat submanifold $\cS$ of the given information manifold $\cQ$, 
there is a unique point $q(\eta)$ in $\cS$ that minimizes 
${\rm KL}(p;q(\eta))$ and it is given by the point where the m-geodesic
from $p$ meets $\cS$ orthogonally with respect to the Fisher information
metric (see Theorem 1 of \cite{IkTaAma}).

\smallskip

Thus, from the point of view of information geometry, it seems especially
interesting to look at cases of the thermodynamic semiring structures
$$ x(\eta) \oplus_{{\rm KL}_{q(\eta)}} y(\eta) =\sum_p \rho^{-{\rm KL}(p;q(\eta))} x(\eta)^p y(\eta)^{1-p} $$
for distributions $q(\eta)$ that vary along e-flat submanifolds of information manifolds
and recast some Riemannian aspects of information geometry in terms of
algebraic properties of the thermodynamic semirings.

\subsection{Tropical geometry}\label{tropicalsec}

Most of our results have a very natural thermodynamic interpretation
when written explicitly in the case of the tropical semifield (seen as
a prototype example of characteristic one semiring as in \cite{Co}, 
\cite{CC}). Thus, besides the original motivation arising in the
context of $\F_1$ geometry, it is possible that the theory of thermodynamic
semirings we developed here may have some interesting applications in
the setting of tropical geometry \cite{ItMiShu}.

\smallskip

The use of tropical geometry in the context of probabilistic inference in
statistical models was recently advocated in \cite{PaStu}. In that approach
one considers polynomial maps from a space of parameters to the space of
joint probability distributions on a set of random variables. These give
statistical models described by algebraic varieties. The tropicalization
of the resulting algebraic variety is then used as a model for parametric
inference, for instance, by interpreting marginal probabilities
as coordinates of points on the variety.

\smallskip

It would therefore seem interesting to extend the encoding of thermodynamic
and information-theoretic properties into the additive structure of the semiring 
to the broader context of tropical varieties. In particular one can consider 
the patchworking process, where operations are peformed on the "quantized" varieties, 
and then the limit in the Maslov dequantization, corresponding to the residue morphism $T \rightarrow 0$, 
is performed, obtaining the new tropical variety.

\smallskip

Observe, for instance, that in the usual setting of tropical geometry, in passing from
an algebraic variety to its tropicalization, starting with a
polynomial $f$ defining a hypersurface $V$ in $(\C^*)^n$,
one can proceed by first considering an associated Maslov dequantization, given by a one-parameter family $f_h$, whose zero set one denotes by $V_h$. One
then considers the amoeba obtained by mappint $V_h$ to $\R^n$ under the
map ${\rm Log}_h(z_1,\ldots,z_n)=(h \log |z_1|, \ldots, h \log |z_n|)$. One 
obtains in this way the amoeba $\cA_h={\rm Log}_h(V_h)$. As we send the
parameter $h\to 0$, the subsets $\cA_h\subset \R^n$ converge in the
Hausdorff metric to the tropical variety ${\rm Tro}(V)$, see \cite{Litv}.
For example, for a polynomial of the form
$f(x)=\sum_k a_k x^k$,  one obtains $f_h(x)$ by passing to $a_k=e^{b_k}$
and $x^k= e^{kt}$, so that one can then replace $v=\log (\sum_k e^{kt+b_k} )$
by the deformed $v_h =h \log (\sum_k e^{(kt+b_k)/h} )$, which in turn defines
the dequantized family $f_h(x)$.

\smallskip

By comparing with Proposition \ref{propShsum}, one can see
that the Maslov dequantization can be expressed in terms of
the operation $\oplus_{{\rm Sh},T}$, where the dequantization
parameter $h$ plays the role of the temperature $T$,
as also observed in \cite{Co}.  Therefore, one can introduce
variants of the Maslov dequantization procedure, based on
other operations $\oplus_{S,T}$, for other choice of information
measures.  In particular,  one
can consider dequantizations based on various
$n$-ary information measures of the form
$$ (x_1 \oplus_{S,T} \cdots \oplus_{S,T} x_n)_{{\bf T}} = 
\min_{\sum_i p_i=1} ( \sum_i p_i x_i - T S_{\bf T} (p_1, \ldots, p_n)), $$
with the data labelled by trees ${\bf T}$, as we 
described in \S  \ref{operadsec} above.

\smallskip

While one can expect that the tropical
limit itself will be independent of the use of different information
measures in the dequantization procedure, operations performed
at the level of the amoebas $\cA_h$ will likely have
variants with different properties when the Shannon entropy
is replaced by other information measures of the kind considered 
in this paper.

\medskip
\subsection{The thermodynamics of $\R^{un}$}

In the characteristic $p$ case, the functoriality of the Witt construction 
provides a way to construct extensions of the field of $p$-adic
numers $\Q_p ={\rm Frac}(\Z_p)$ using the fact that $\Z_p = \cW_p(\F_p)$,
and applying the same Witt functor to extensions $\F_q$. This gives
$\cW_p(\F_q)= \Z_p[\zeta_{q-1}]$, which is the valuation ring of an  
unramified extension $\Q_p(\zeta_{q-1})$ of $\Q_p$, see \cite{Lor}.

\smallskip

It was observed in \S 7 of \cite{Co} that, in the case of the
chracteristic one version of the Witt construction, when one
considers the $\oplus_{{\rm Sh},T}$ simultaneously for all possible 
temperatures $T$, one can describe a candidate analog of
``unramified extension" $\R^{un}$ in terms of analogs of
Teichm\"uller characters given in the form $\chi_T(f)=f(T)^{1/T}$
and an analog of the residue morphism of the form 
$\epsilon(f)=\lim_{T\to 0} \chi_T(f)(T)^T$.
\smallskip

We can formulate this in the general case. We find, first of all that the Frobenius lifts do not depend on the information measure.

\begin{prop}
If $R$ is a thermodynamic semiring over a suitably nice semifield $K$ defined by the information measure $S$, 
then the Frobenius lifts from $K$ to $R$ in such a way that
\[F_r(x(T)) = x(T/r)^r. \]
\end{prop}
\proof
We see that this is a result of the general form of the temperature dependence in the current context. In symbols, we are looking for
\[F_r(x(T)\oplus_S y(T)) = \sum e^{r f(T) S(\alpha)}x(f(T))^{r \alpha}y(f(T))^{r (1-\alpha)} \]\[ = \sum e^{f(T) S(\alpha)}x(f(T))^{r \alpha}y(f(T))^{r (1-\alpha)}=F_r(x(T))\oplus_S F_r(y(T)),\]
where the residue morphism forces $F_r(x(T)) = x(f(T))^r$ for some invertible $f$, depending on $r$. We see from the above that $f(T) = T/r$, proving the claim.
\endproof

This forces the characters to have the same form as in the Shannon entropy case, ie. $\chi_T (f)(T) = f(T)^{1/T}$. However, these characters are additive only if $(x(T) \oplus_S y(T))^{1/T} = x^{1/T}+y^{1/T}$, which means $S$ must produce the same 
thermodynamic structure as the Shannon entropy, hence, by a theorem above, $S$ is the Shannon entropy. Note that this analysis holds also in the $q$-deformed Witt construction leading to the Tsallis entropy discussed in \S \ref{WittTsallisSec}.

\smallskip

If we pass to the field of fractions of these characters, and consider further infinite sums of these characters, 
the resulting expressions begin to resemble partition functions in the Euclidean path integral formulation, see \S 7 in \cite{Co}. Indeed, if one uses instead $\bR^{{\rm min}, +} \cup \{\infty \}$, these are equal to equilibrium free energies of the type observed in \S \ref{thermoSec}. The failure of the additivity of the characters in $\cR^{\rm un}$ can thus be interpreted in terms of nonextensivity. This suggests that, as this candidate for $\cR^{\rm un}$ is investigated, new algebraic interpretations of nonextensivity will arise. It would also be interesting to see if a notion of character which is additive on the $q$-deformed Witt construction could give rise to a one-parameter family of $\cR^{\rm un}$'s.

\subsection{Thermodynamics in positive characteristics}

The main motivation for the Witt construction in characteristic one given in 
\cite{Co} and \cite{CC}, which provides the prototype example of a thermodynamic 
semiring built on the Shannon entropy, is to provide an analog in characteristic one
of the formulae for the summation of Teichm\"uller representatives in the case of 
multiplicative lifts to $\Z_p$ of the characteristic $p$ elements in $\F_p$.

\smallskip

One can then reverse the point of view and start from the more general thermodynamic
semirings associated to other forms of entropy, such as R\'enyi, Tsallis, Kullback--Leibler,
with their axiomatic characterizations, and look for characteristic $p$ analogs of
non-extensive thermodynamics and other such variants of statistical physics.

\smallskip

For instance, we saw in \S \ref{tsallisSec} above that there is a one-parameter
deformation of the Witt construction in characteristic one, which yields a 
characterization of the Tsallis entropy ${\rm Ts}_\alpha$ as the unique binary
information measure that satisfies the associativity, commutativity and unity
constraints for this deformed $\oplus_{S,T,\alpha}$ operation. 

\smallskip

One thinks of the original $\oplus_{{\rm Sh},T}$ with the Shannon
entropy as in \cite{Co} and \cite{CC}, as being the correct analog
in characteristic one of the $p$-adic Witt construction 
$$ x\oplus_w y = \sum_{s\in I_p} w_p(s) x^s y^{1-s}, $$
with $I_p$ the set of rational numbers in $[0,1]$ with denominator a 
power of $p$ and 
$$ w_p(s)= \sum_{a/p^n =s} w(p^n,a) T^n \in \F_p((T)), $$
where the $w(p^n,k) \in \Z/p\Z$, for $0<k<p^n$ are determined by
the addition of Teichm\"uller representatives
$$ \tau(x)+\tau(y) =\tau(x+y) 
+\sum_{n=1}^\infty \tau\left( \sum w(p^n,k) x^{k/p^n} y^{1-k/p^n}\right) p^n . $$

\smallskip

Thus, one can equivalently think of the universal sequence of the $w(p^n,k)$ as being
the characteristic $p$ analog of the Shannon information. Adopting this viewpoint,
one would then expect that the one-parameter deformation of the Witt construction
in characteristic one described in \S \ref{tsallisSec}, which leads naturally from the
Shannon entropy to the non-extensive Tsallis entropy, may correspond to
an analogous deformation of the original $p$-adic Witt construction that
leads to a notion of non-extensive entropy and non-ergodic thermodynamics
in characteristic $p$. 

\smallskip

It should be mentioned that there are in fact interesting known $q$-deformations 
of the Witt constructions, see for instance \cite{Oh}. These can naturally be 
described within the setting of $\Lambda$-rings (see \cite{Oh}). This seems
especially useful, in view of the whole approach to $\F_1$ geometry based
on $\Lambda$-rings, developed by James Borger in \cite{Borg} and
\cite{Borg2}, \cite{Borg3} (see also \cite{Man}, \cite{Mar} for other related viewpoints). However, a reader familiar with the positive characteristic Witt construction will notice that Connes and Consani's construction generalize the $p$-Witt ring from a rather unconventional expression for its addition. This is difficult to translate into the $\Lambda$-ring approach to the Witt ring. A definition of $\Lambda$-rings in characteristic one which reproduces the Witt rings considered in this paper would likely be interesting both geometrically and physically.
\smallskip

This also suggests that identifying suitable analogs of other entropy
functions (Tsallis, R\'enyi, Kullback--Leibler) in characteristic $p$,
via deformations of the ring of Witt vectors, may also further our 
understanding of $\F_1$-geometry in the $\Lambda$-ring approach.

\bigskip
\bigskip

{\bf Acknowledgment.} This paper is based on the results of  the second author's
summer research project, supported by the Summer Undergraduate Research
Fellowship program at Caltech. The first author is partly supported by 
NSF grants DMS-0901221 and DMS-1007207.

\end{document}